\documentclass[11pt]{article}
\usepackage{amsmath}
\usepackage{amssymb}
\usepackage{amsthm}
\usepackage{amsfonts}
\usepackage{algorithm}
\usepackage{algorithmic}
\usepackage{graphicx}
\usepackage{pdfpages}
\usepackage{xcolor,colortbl}
\usepackage{multicol}
\usepackage{url}
\usepackage{subfigure}
\usepackage{ulem}
\usepackage{listings}
\usepackage{tikz}
\definecolor{Gray}{gray}{0.92}
\definecolor{LightCyan}{rgb}{0.88,1,1}

\newcolumntype{a}{>{\columncolor{Gray}}c}
\newcolumntype{b}{>{\columncolor{white}}c}

\usepackage{bbm,epsfig,graphics,epic,color,rotating,color}
\numberwithin{equation}{section}

\def\E{{\mathbb E}}
\def\P{{\mathbb P}}

\newtheorem*{theorem*}{Theorem}
\newtheorem{theorem}{Theorem}[section]

\newtheorem{assum}[theorem]{Assumptions}

\theoremstyle{definition}

\makeatletter
\def\BState{\State\hskip-\ALG@thistlm}
\makeatother

\title{The Monty Hall problem revisited}
\date{Departamento de Estad\'i{}stica\\ Universidad del B\'io-B\'io. }
\author{Christian Caama\~no-Carrillo\footnote{Avda. Collao 1202, CP 4051381, Concepci\'on, Chile. E-mail address: chcaaman@ubiobio.cl} \ and Manuel Gonz\'alez-Navarrete\footnote{Corresponding author. E-mail address: magonzalez@ubiobio.cl}}

\begin{document}

\maketitle %
\thispagestyle{empty} %
\baselineskip=14pt

\vspace{2pt}

\begin{abstract}
We propose a new approach to solve the classical Monty Hall problem in its general form. The solution is based on basic tools of probability theory, by defining three elementary events which decompose the sample space into a partition. The probabilities of each element of the partition allow us to compute the conditional and marginal probabilities of winning.
\end{abstract}


%
\bigskip
%


%
\section{Introduction}
\label{sec:intro}
The Monty Hall problem appeared first in a letter by S. Selvin to the American Statistician \cite{se1,se2}, it is a nice and controversial problem for introductory courses in probability, statistics, and game theory. The problem can be posed as follows:
\medskip

\begin{minipage}[c]{16cm}
\textit{You are playing a game on a TV show. There are three doors. One of them has a car in
the backside, and the other ones have goats. You select one of the doors, say door No. 1. Before opening it, the host (who knows where the car is) opens another door, say No. 3, which has a goat. And then he gives you the possibility to change, allowing you to pick door No. 2. What do you do?}
\end{minipage}
\medskip

Nowadays, the answer for this problem is well known: the best strategy is changing the selected door. However, in the early 90's a publication of M. Vos Savant \cite{vs1} generated a controversy, because of the answer (solution) given to the problem. Essentially, the article exposed the intuitive arguments of the right (mathematical) solution proposed by Selvin \cite{se1}. Nevertheless, some mathematicians did not agree with Vos Savant' explanation and stated that it does not matter if you change the selected door. Since there are two remaining doors, the probability of door No. 1 (or No. 2) having a car is $1/2$. Subsequent articles \cite{vs2,vs3} discussed some opinions that defended the latter argument, these criticized Vos Savant even saying something like:
\medskip

\begin{minipage}[c]{16cm}
\textit{May I suggest that you obtain and refer to a standard textbook on probability before you try to answer a question of this type again?}
\end{minipage}
\medskip

In fact, a standard textbook on probability theory \cite{BT,Fe,LR} is enough to confirm Vos Savant' solution. Clearly, she was sure about that and in \cite{vs2} encouraged the readers to check empirically her answer to the Monty Hall problem. Even more, Vos Savant \cite{vs1} suggests 1.000.000 doors rather than 3. In this case, there are 999.999 doors with goats behind them and one door with a prize. After the player picks a door, the host opens 999.998 of the remaining doors. On average, in 999.999 times out of 1.000.000, the remaining door will contain the prize.

In this sense, the literature presents several solutions \cite{se2,MC,cy,GiL,car,gsl,luc,Gi2}.
We remark that, the Monty Hall problem is relevant in the teaching of probability, because of the misperception that it can generate. In this sense, initial concepts from theory of probability can be studied by the statement of this problem. For instance, conditional probability and the law of total probability. However, a set theory approach is also of interest, because of the several possibilities to define the events of interest. Our proposal is based in the definition of the following three events:

\medskip

$E$: the contestant initially chooses  the correct door.

$C$: the contestant changes the selected door.

$W$: the contestant wins the car.

\medskip

\begin{figure}[!hbtp]
\begin{center}
\includegraphics[width=0.2\linewidth, height=0.15\textheight]{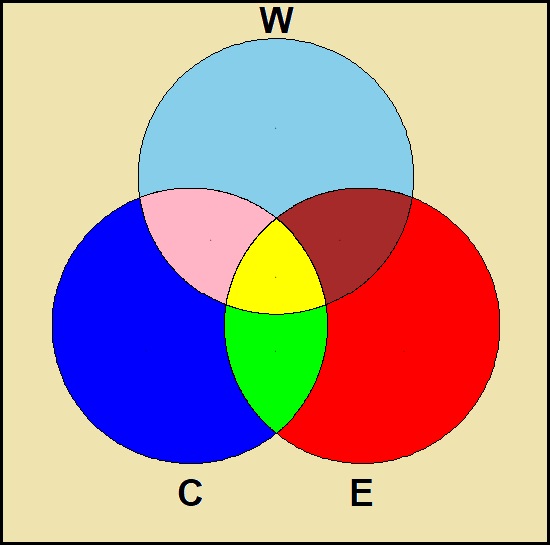}
\end{center}
	\caption{Venn diagram associated to the three events above.}
	\label{vendiamgramevent}
\end{figure}

At this point the notion of partition of a set is relevant to understand the solution that we propose (see for details \cite{BN}). A partition divides a set into separate and non-empty events. For instance, Figure \ref{vendiamgramevent} includes a Venn diagram representing a partition with eight events, where the union of these subsets is equal to the sample space. The same collection is shown in the tree diagram in Figure \ref{tree}.

\tikzstyle{level 1}=[level distance=3cm, sibling distance=4.5cm]
\tikzstyle{level 2}=[level distance=3cm, sibling distance=2.3cm]
\tikzstyle{level 3}=[level distance=3cm, sibling distance=1.2cm]

\tikzstyle{bag} = [text width=4em, text centered]
\tikzstyle{end} = [minimum width=3pt, inner sep=0pt]

\begin{figure}
  \centering

\begin{tikzpicture}[grow=right]
  \node[bag] {}
    child {
        node[bag] {\includegraphics[width=1cm]{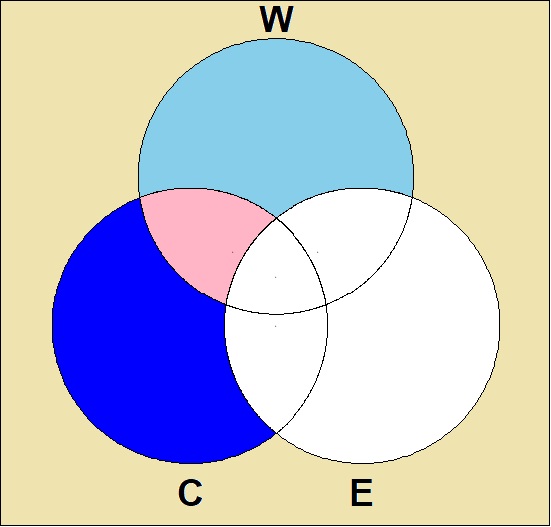}}
            child {
             node[bag] {\includegraphics[width=1cm]{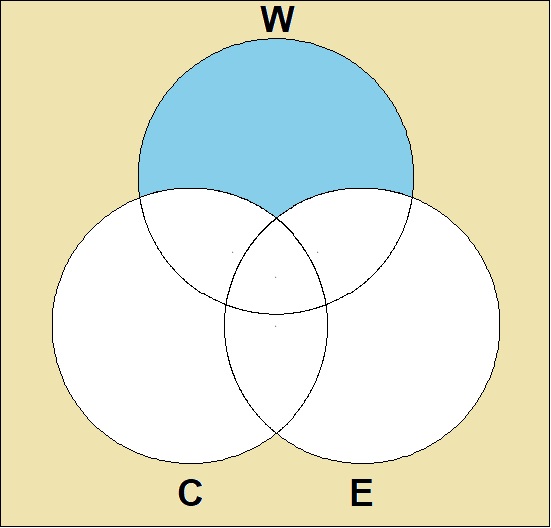}}
              child {
                node[end, label=right:
                   {\includegraphics[width=1cm]{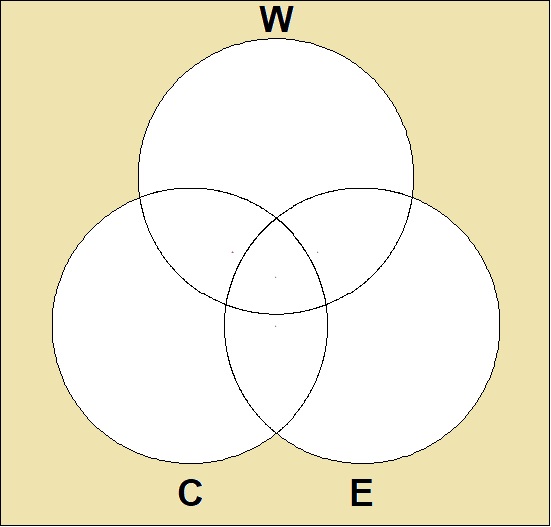}}] {}
                edge from parent
                node[above] {$W^{c}$}
                node[below]  {}
            }
            child {
                node[end, label=right:
                    {\includegraphics[width=1cm]{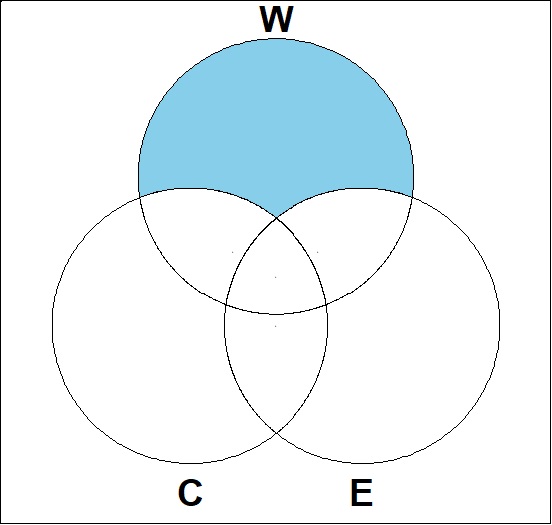}}] {}
                edge from parent
                node[above] {$W$}
                node[below]  {}
            }
                edge from parent
                node[above] {$C^{c}$}
                node[below]  {}
            }
            child {
             node[bag] {\includegraphics[width=1cm]{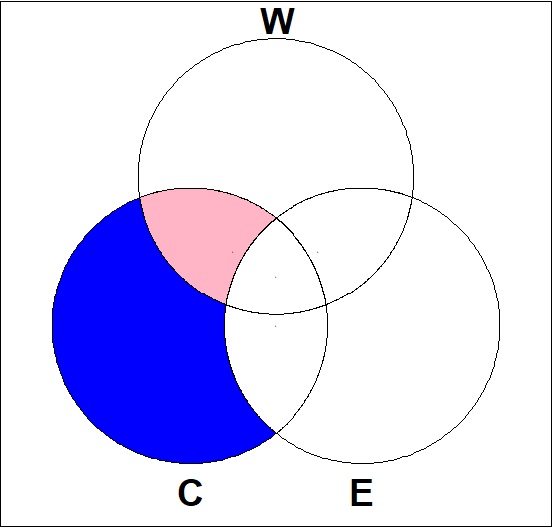}}
              child {
                node[end, label=right:
                   {\includegraphics[width=1cm]{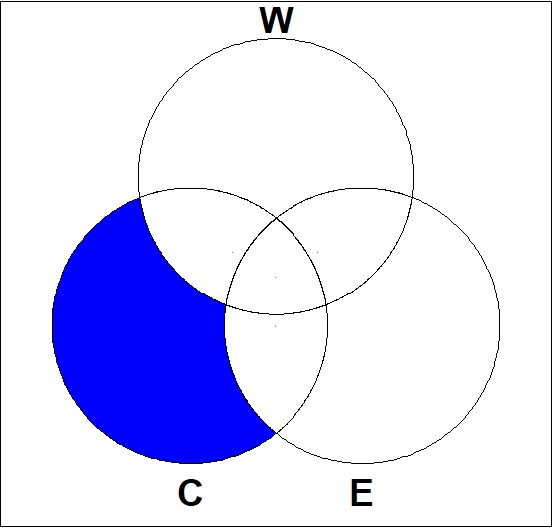}}] {}
                edge from parent
                node[above] {$W^{c}$}
                node[below]  {}
            }
            child {
                node[end, label=right:
                    {\includegraphics[width=1cm]{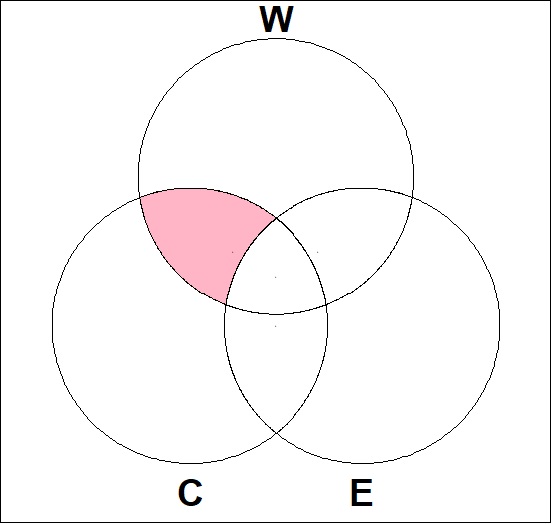}}] {}
                edge from parent
                node[above] {$W$}
                node[below]  {}
            }
                edge from parent
                node[above] {$C$}
                node[below]  {}
            }
            edge from parent
            node[above] {$E^{c}$}
             node[below]  {}
    }
    child {
        node[bag] {\includegraphics[width=1cm]{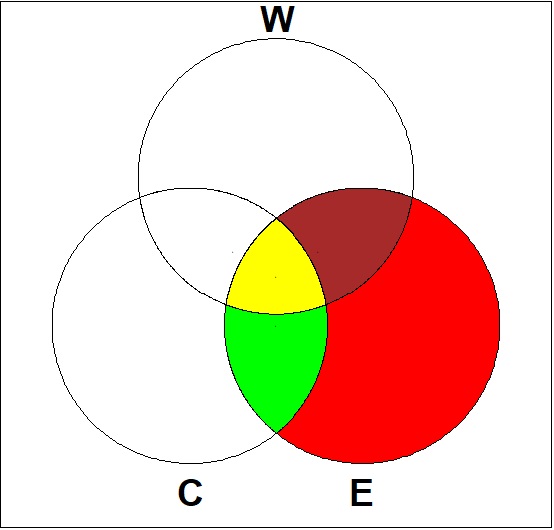}}
            child {
             node[bag] {\includegraphics[width=1cm]{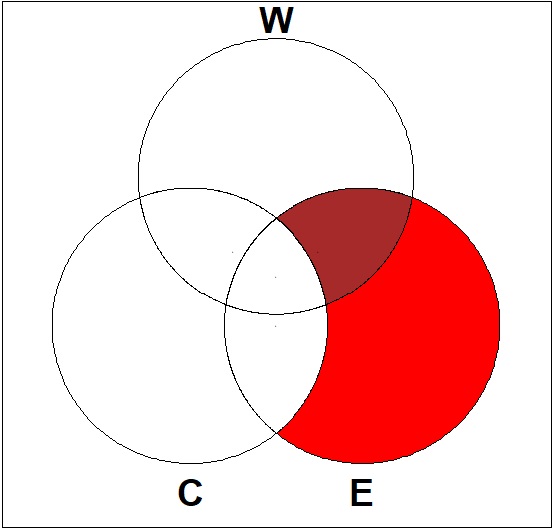}}
              child {
                node[end, label=right:
                   {\includegraphics[width=1cm]{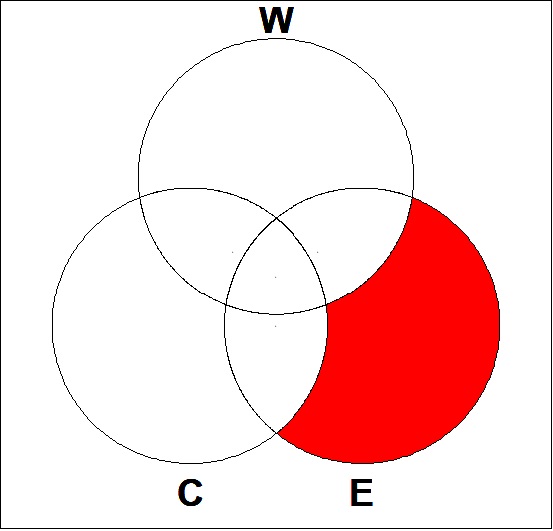}}] {}
                edge from parent
                node[above] {$W^{c}$}
                node[below]  {}
            }
            child {
                node[end, label=right:
                    {\includegraphics[width=1cm]{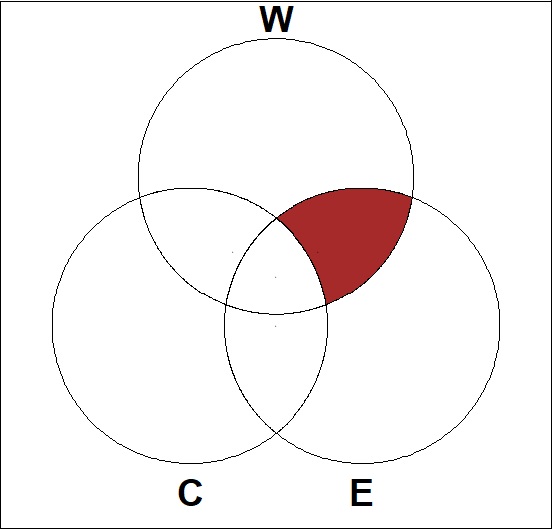}}] {}
                edge from parent
                node[above] {$W$}
                node[below]  {}
            }
                edge from parent
                node[above] {$C^{c}$}
                node[below]  {}
            }
            child {
             node[bag] {\includegraphics[width=1cm]{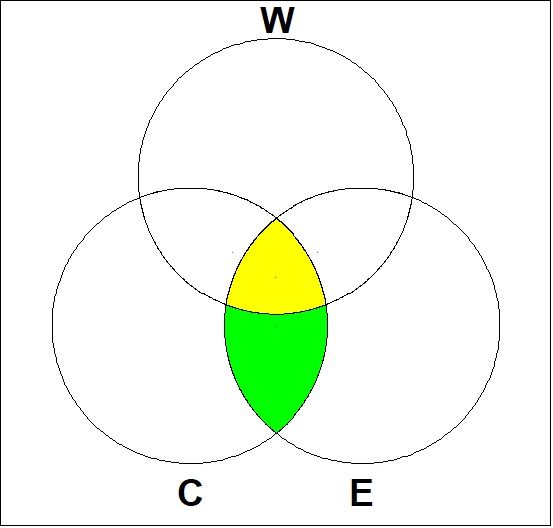}}
              child {
                node[end, label=right:
                   {\includegraphics[width=1cm]{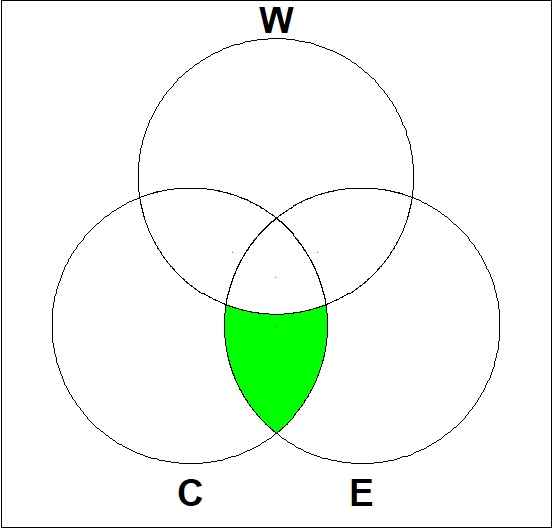}}] {}
                edge from parent
                node[above] {$W^{c}$}
                node[below]  {}
            }
            child {
                node[end, label=right:
                    {\includegraphics[width=1cm]{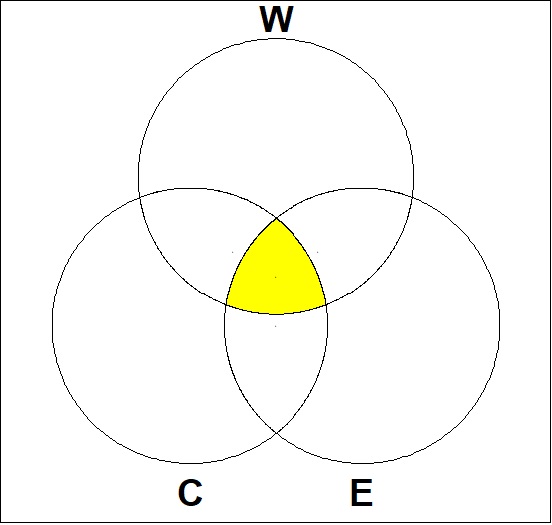}}] {}
                edge from parent
                node[above] {$W$}
                node[below]  {}
            }
                edge from parent
                node[above] {$C$}
                node[below]  {}
            }
            edge from parent
            node[above] {$E$}
             node[below]  {}
    };
\end{tikzpicture}
  \caption{Tree diagram representing the events and their relation with the Venn diagram in Figure \ref{vendiamgramevent}.}\label{tree}
\end{figure}

In this paper we discuss the implicit assumptions made by the different approaches in the literature and discussed their importance in the solution of the problem (see \cite{MGN}). The advantage of our proposal is its simplicity. The definition of the three above events  allow us to solve the problem with elementary probabilistic tools, turning the problem a candidate to introduce the concepts of conditional probability in a first undergraduate course or even more, in a lecture for secondary school.

The theoretical solutions are then checked by numerical simulations and the pseudo-codes are included to allow the replication of our results. The goal is to mimic the empirical evidence supporting the solutions, as required by Vos Savant in \cite{vs1}.

The rest of the paper is organized as follows. In Section \ref{sec:model} we include the formulation of the probability models and present their theoretical solutions; Section \ref{sec:simu} is dedicated to the numerical simulations of the solutions, including the arguments to choose the numbers of iterations. Finally, Section \ref{sec:comment} concludes providing final remarks.


\section{The probability model and its solution: computing conditional and marginal probabilities of winning}
\label{sec:model}

We propose an approach which allows us to solve the MH problem in the original form and its generalizations as proposed in \cite{se2} and \cite{vs2}, that is, the existence of $n$ doors, with $n \ge 3$.

    We use the chain rule in probability to obtain, for instance

    \begin{equation}
    \label{chain}
\P(E\cap C \cap W) =  \P ( W| E\cap C) \cdot \P(C| E)\cdot \P(E),
\end{equation}
also for $\{E\cap C \cap W^c\}, \{E\cap C^c \cap W\}$ and so on, for all eight subsets in the partition proposed by Venn diagram or for each of the branches in the tree diagram in Figure \ref{tree}.

\begin{assum}
\label{assum}
We solve the MH problem assuming that initially the contestant selects the door uniformly at random, that is, $\P(E) = \frac{1}{n}$, with $n \ge 3$. Moreover, we suppose that the contestant does not know if the initially chosen door has the car and takes the decision to change the door with a given probability $p \in [0,1]$, then, $\P(C| E) = \P(C| E^c)=p$.
\end{assum}

Then, we compare the conditional probabilities

    \begin{equation}
 \P ( W| C) \text{ and }   \P ( W| C^c),
\end{equation}
of winning the prize by changing and winning by staying, respectively and concluding which strategy is better. In addition, we obtain the marginal probability of winning, $\P(W)$.

This assumption generalizes previous solutions which implicitly assumed equal probabilities, $p=1/2$. That is, the contestant has not a preference to change or not, making his (her) decision at random. For details, see comments in Section \ref{sec:comment}.


We solve the MH problem and generalize it in the sense proposed by Vos Savant \cite{vs1}. Then we suppose there are $n$ doors and consider two different strategies assumed by the host to open door(s) after the contestant initial choice. In the first case, in Section \ref{sec:MH2}, the host opens $(n-2)$ doors, then the contestant has to choose from the remaining 2 doors. In the second case, Section \ref{sec:MH3}, the host opens just one door and then the contestant has to choose from the remaning $(n-1)$ doors. In the original MH problem, with 3 doors, both strategies coincide.

In Section \ref{sec:simu} we will develop numerical simulations of both strategies, aiming to check the theoretical solution, as suggested by Vos Savant to her readers in \cite{vs1}.


\subsection{Original MH problem}
\label{sec:MH1}

There are 3 doors, then by Assumptions \ref{assum},
   \begin{equation}
\P(E) =  \frac{1}{3} \text{ and } \P(C | E) =  p.
\end{equation}

Therefore, the probability in \eqref{chain} is computed by

    \begin{equation}
\P(E\cap C \cap W) = 0 \cdot  p \cdot  \frac{1}{3},
\end{equation}
given that if the contestant chooses the correct door and changes it, he is finally selecting a wrong door, that is, losing with probability 1. Now, we obtain
    \begin{equation}
\P(E\cap C \cap W^c) =  \P ( W^c| E\cap C) \cdot \P(C| E)\cdot \P(E) = 1 \cdot p \cdot  \frac{1}{3}.
\end{equation}
Moreover,

\begin{equation}
\begin{array}{c}
\P(E^c\cap C^c \cap W^c)= \P ( W^c| E^c\cap C^c) \cdot \P(C^c| E^c)\cdot \P(E^c)  = \frac{2}{3}(1-p),\\[0.3cm]
\P(E\cap C^c \cap W) =\P ( W| E\cap C^c) \cdot \P(C^c| E)\cdot \P(E) =\frac{1}{3}(1-p),\\[0.3cm]
\P(E^c\cap C \cap W)=\P ( W| E^c\cap C) \cdot \P(C| E^c)\cdot \P(E^c) =\frac{2}{3} p,
\end{array}
\end{equation}
and by similar arguments
    \begin{equation}
\P(E\cap C^c \cap W^c)= \P(E^c\cap C^c \cap W)=\P(E^c\cap C \cap W^c)=0.
\end{equation}

Therefore, the conditional probability of winning by changing the door is
    \begin{equation}
    \label{pc}
\P(W|C) = \frac{\P(W \cap C)}{\P(C)} = \frac{\P(E\cap C \cap W) + \P(E^c\cap C \cap W)}{\P(C)} = \frac{0 + \frac{2}{3} p }{p} = \frac{2}{3},
\end{equation}
and the conditional probability of winning by keeping the door is

    \begin{equation}
    \label{pcc}
\P(W|C^c) = \frac{\P(W \cap C^c)}{\P(C^c)} = \frac{\P(E\cap C^c \cap W) + \P(E^c\cap C^c \cap W)}{\P(C^c)} = \frac{ \frac{1}{3}(1- p) + 0 }{1-p} = \frac{1}{3}.
\end{equation}

In this sense, the probabilities in \eqref{pc} and \eqref{pcc} are in acccordance with the solution given by Selvin \cite{se1} and Vos Savant \cite{vs1}. That is, the best strategy is changing the initially selected door.

Finally, by using the law of total probability

\begin{equation}
\label{pw1}
\begin{array}{lll}
\P(W) &= & \P ( W| C) \P(C) + \P ( W| C^c) \P(C^c)\\[0.3cm]
& =  & \displaystyle\frac{1}{3}  + \displaystyle\frac{1}{3} p.
\end{array}
\end{equation}

This probability of winning is a linear function of the parameter $p$ (probability of change). In particular, the intercept and the slope of this function are the same ($1/3$). This fact implies that as $p$ grows the probability in \eqref{pw1} increases. See Figure \ref{fig:f1}.


\subsection{MH with $n$ doors, the host leaves two closed doors}
\label{sec:MH2}

There are $n$ doors and the assumption is kept.

   \begin{equation}
\P(E) =  \frac{1}{n} \text{ and } \P(C | E) =  p.
\end{equation}

In this case, all conditional probabilities $\P(W| E\cap C),\P(W| E\cap C^c),\P(W| E^c\cap C)$ and $\P(W| E^c\cap C^c)$ remain the same as in previous case. Finally, the conditional probabilities of winning by changing or keeping the door are, respectively

    \begin{equation}
\P(W|C) =   \frac{n-1}{n} \text{ and } \P(W|C^c) = \frac{1}{n}.
\end{equation}

Therefore $\P(W|C)  = \P(E^c)$ and $\P(W|C^c)  = \P(E)$, given the dichotomy in the final selection, that is, one door with a car and another with a goat.

Finally

\begin{equation}
\label{pw2}
\P(W) =  \displaystyle\frac{1}{n} + \displaystyle\frac{(n-2)}{n} p .
\end{equation}

\begin{figure}[!hbtp]
\begin{center}
\begin{tabular}{ccc}
\includegraphics[width=5.1cm]{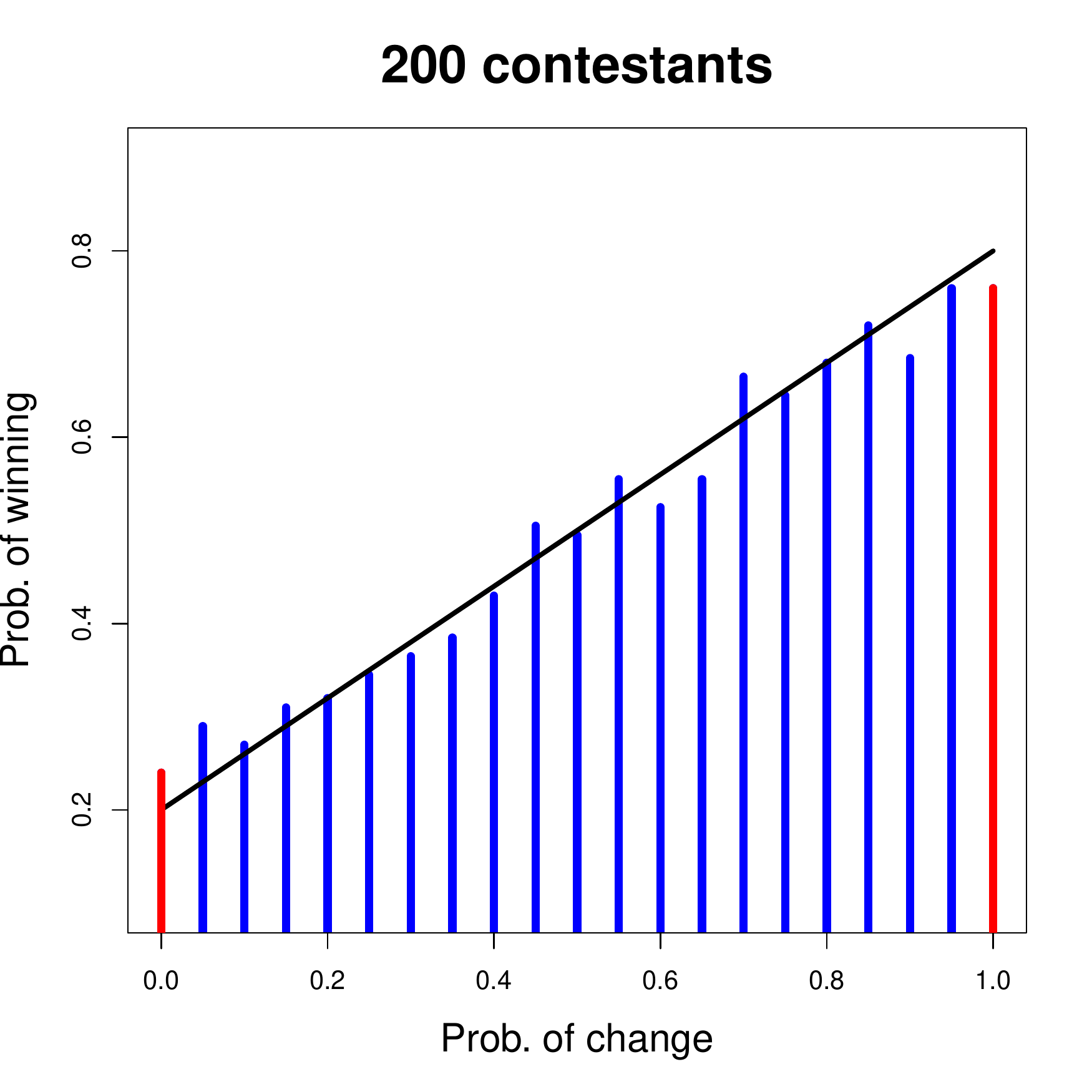} &  \includegraphics[width=5.1cm]{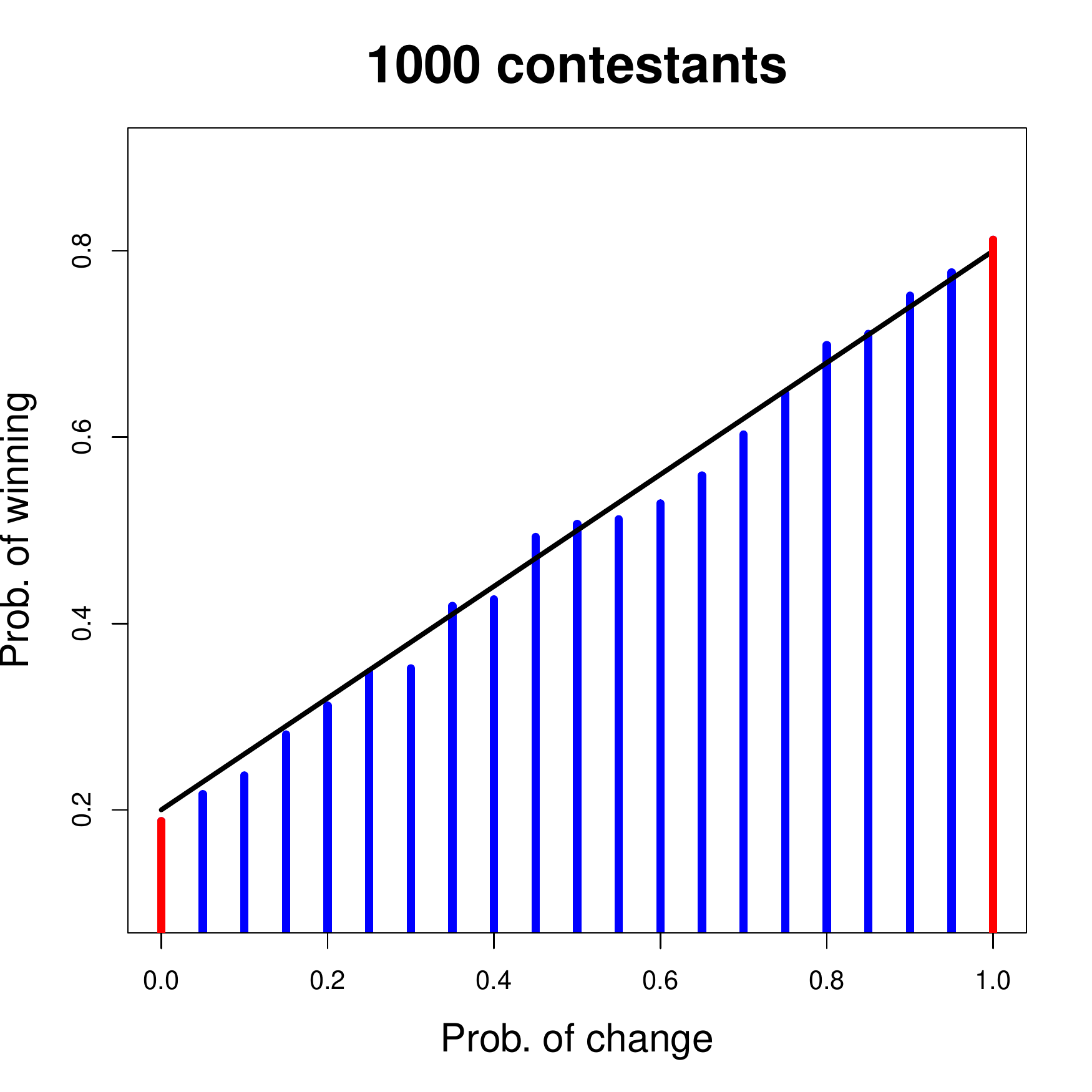} &  \includegraphics[width=5.1cm]{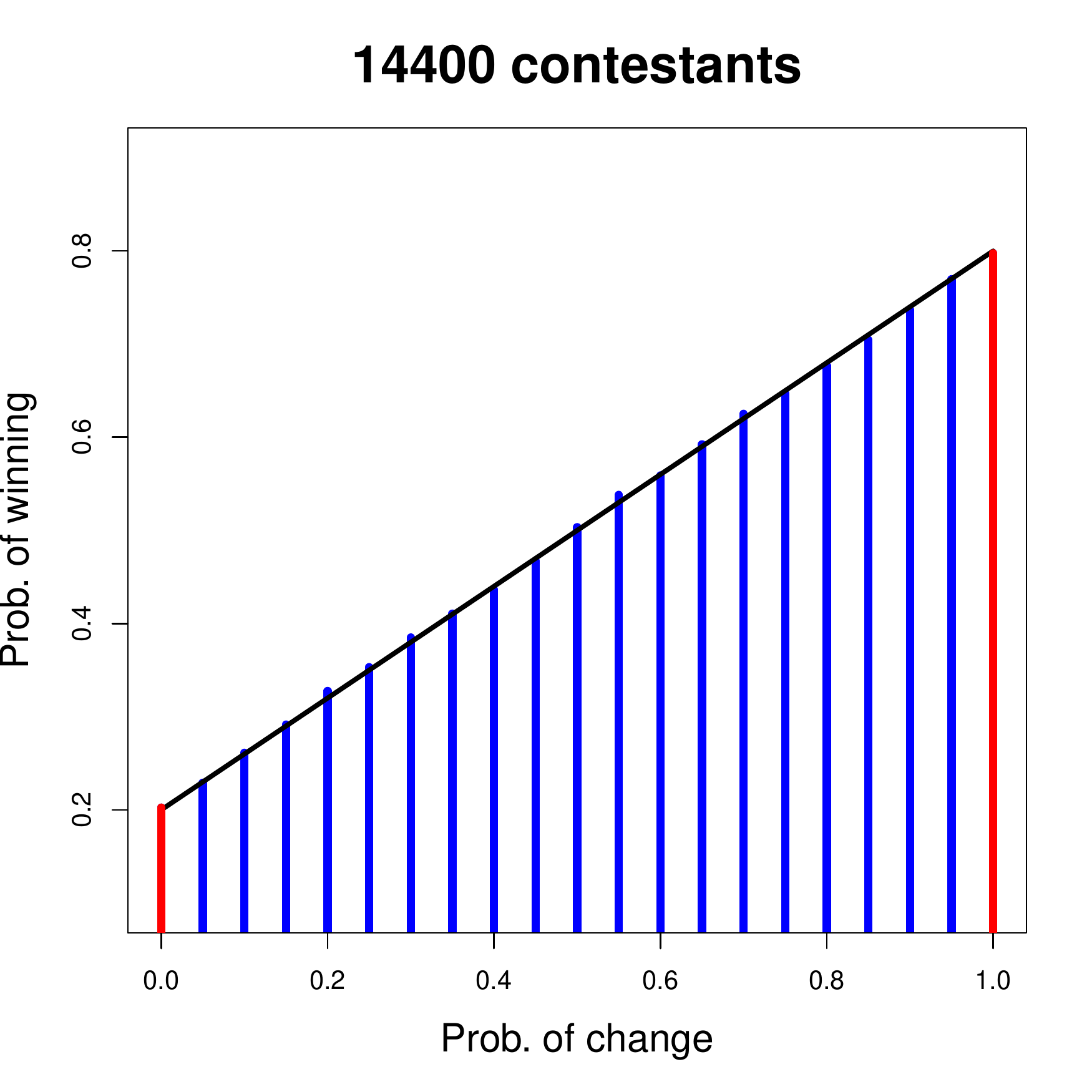}
\\
\includegraphics[width=5.1cm]{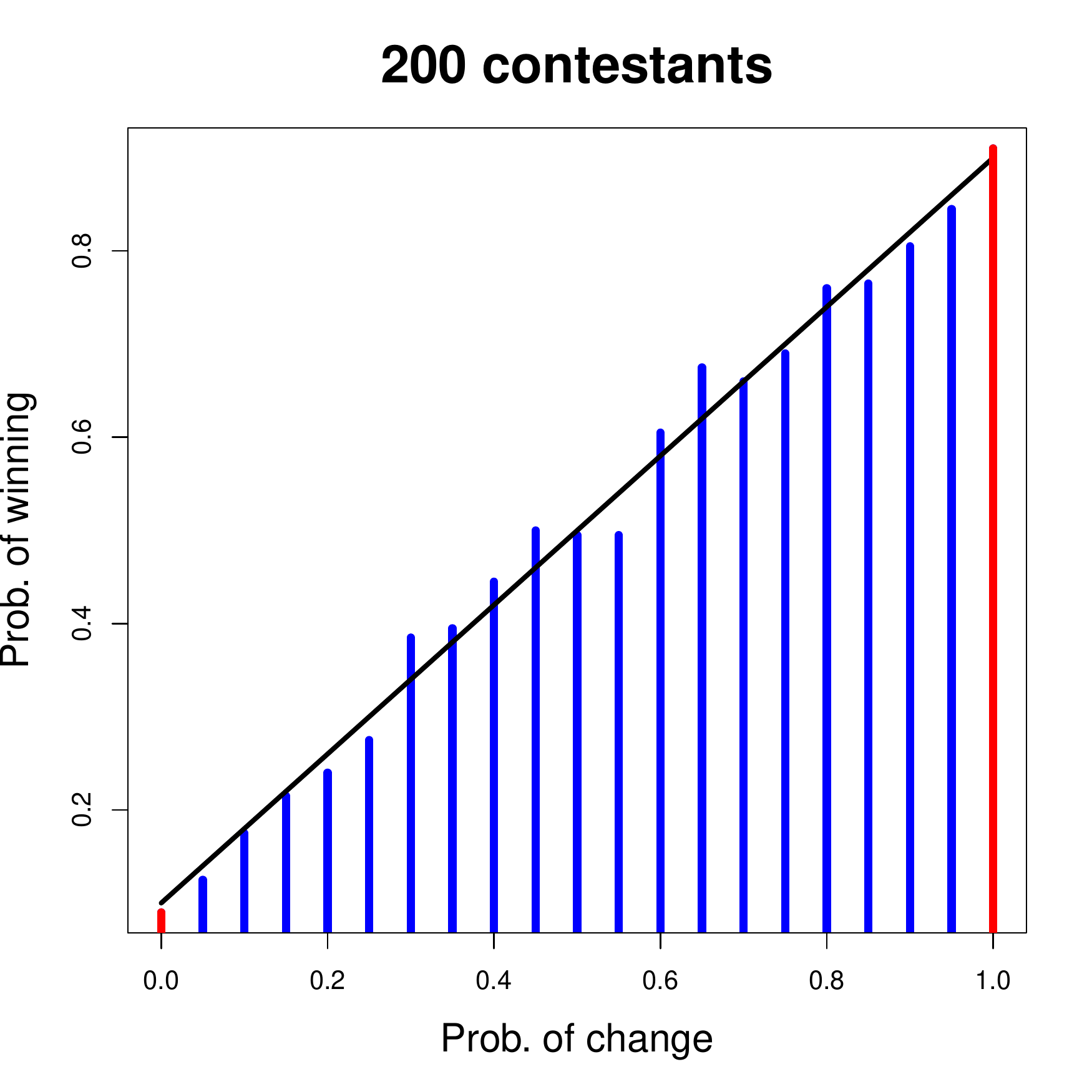} &  \includegraphics[width=5.1cm]{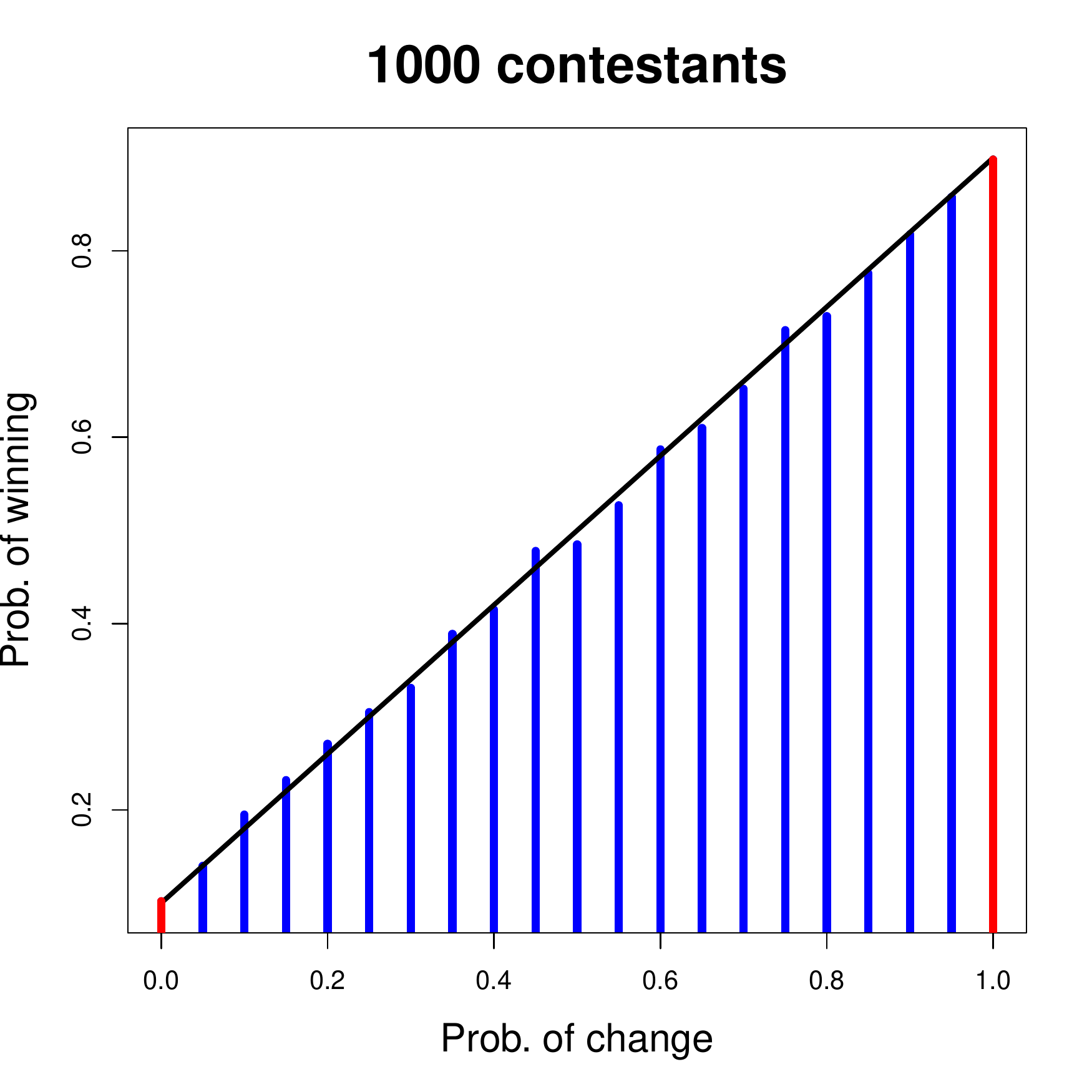} &  \includegraphics[width=5.1cm]{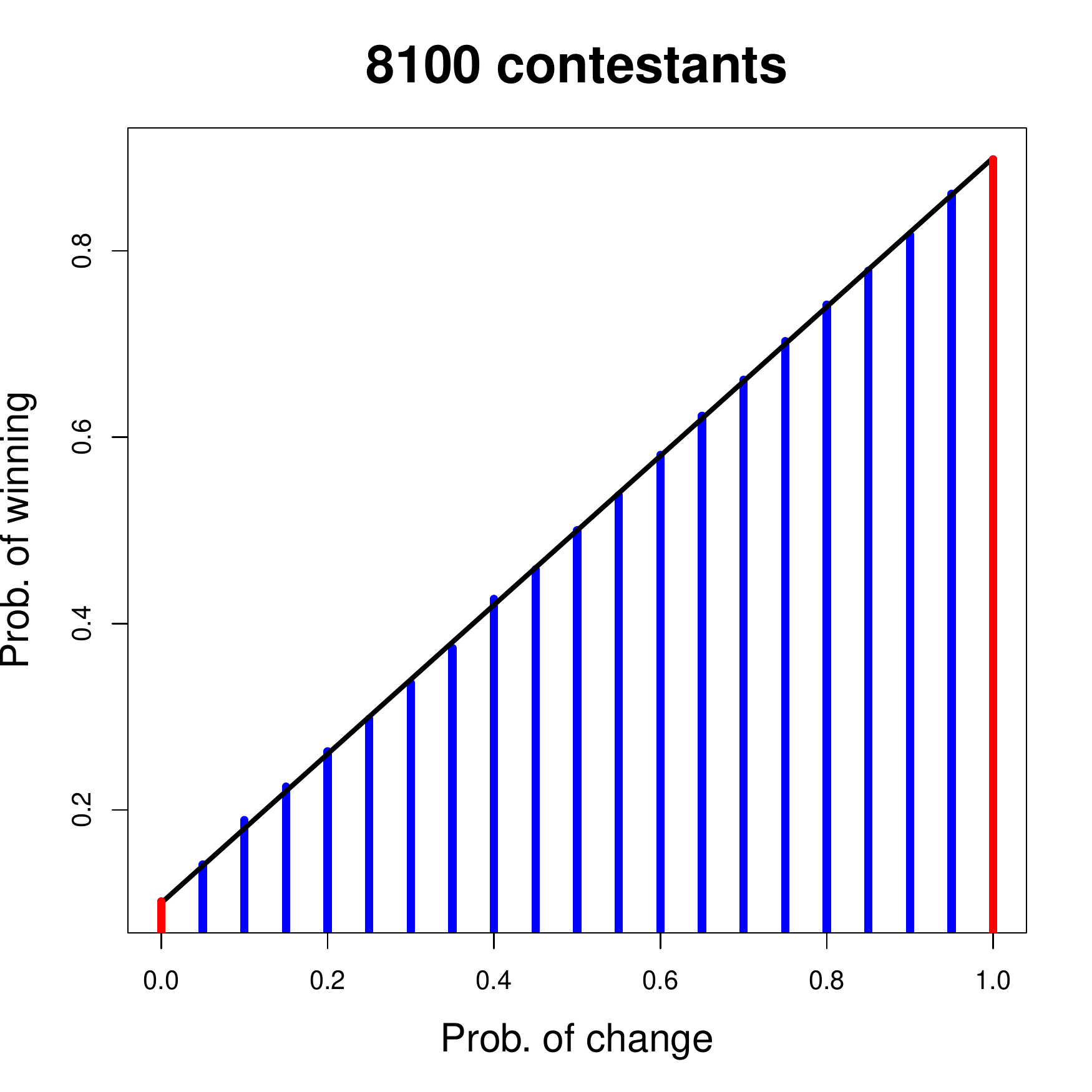}
\end{tabular}
\caption{
Examples of empirical probabilities of winning for the problem in Section \ref{sec:MH2}. In last column we use the optimal $l_0$ given by CLT. Upper part with 5 doors. Bottom part with 10 doors.
 } \label{fig:f2}
\end{center}
\end{figure}

As in previous case, the probability of winning is a linear function of parameter $p$. Of course, for $n=3$ \eqref{pw2} equals to \eqref{pw1}. However, in this situation, both the intercept and the slope depend on the number of doors. In this sense, as $n$ grows to infinity, the intercept goes to zero and the slope is going to one. In other words, for a large number of doors, the probability of winning coincides with the probability of change. That is, $\lim\limits_{n\to \infty} \P(W) = p$. For details see Figure \ref{fig:f2}.

\subsection{MH with $n$ doors, the host opens one door}
\label{sec:MH3}

This problem was mentioned in \cite{se2}, credited to D.L. Ferguson. There are $n$ doors and the host opens one door. In this case

 \begin{equation}
\P(E) =  \frac{1}{n} = 1 - \P(E^c),
\end{equation}
and by Assumptions \ref{assum} $\P(C| E ) = \P(C | E^c) = p $. In addition, we obtain

   \begin{equation}
\P(W| E \cap C) =  \P(W| E^c \cap C^c) =  \P(W^c| E \cap C^c) = 0.
\end{equation}
Which in turn implies that
   \begin{equation}
\P(W^c| E\cap C) =  \P(W^c| E^c\cap C^c) =  \P(W| E\cap C^c) = 1.
\end{equation}
Now, by assuming that the contestant decides to change the door and chooses uniformly at random the new door, we compute
   \begin{equation}
\P(W| E^c\cap C) =  \frac{1}{n-2},
\end{equation}
that is, there is one door with the car belonging to the $n-2$ remaining doors (eliminating the opened and the initially selected doors). Therefore, we also have
   \begin{equation}
\P(W^c| E^c\cap C) = 1 - \P(W| E^c\cap C) = \frac{n-3}{n-2}.
\end{equation}

Finally, we compute the conditional probabilities of winning by changing or keeping the initial door, given by

    \begin{equation}
\P(W|C) = \frac{\P(W \cap C)}{\P(C)} = \frac{\P(E\cap C \cap W) + \P(E^c\cap C \cap W)}{\P(C)} = \frac{0 + \frac{n-1}{n(n-2)} p }{p} = \frac{n-1}{n(n-2)},
\end{equation}
and
    \begin{equation}
\P(W|C^c) = \frac{\P(W \cap C^c)}{\P(C^c)} = \frac{\P(E\cap C^c \cap W) + \P(E^c\cap C^c \cap W)}{\P(C^c)} = \frac{ \frac{1}{n}(1- p) + 0 }{1-p} = \frac{1}{n},
\end{equation}
respectively. In the case of the probability of winning, we obtain

\begin{equation}
\label{pw3}
\begin{array}{lll}
\P(W) &= & \P ( W| C) \P(C) + \P ( W| C^c) \P(C^c)\\[0.3cm]
& =  & p \displaystyle\frac{n-1}{n(n-2)} + (1-p) \displaystyle\frac{1}{n}\\[0.4cm]
& =  &  \displaystyle\frac{1}{n} +  \displaystyle\frac{1}{n(n-2)} p.
\end{array}
\end{equation}
Note that, as in previous cases, this probability is a linear function of parameter $p$. In this sense, as $n$ grows to infinity, the intercept and the slope go to zero. That is, for a large number of doors, the probability of winning goes to zero independently of the contestant decision to change the initially selected door or not. Then, $\lim\limits_{n\to \infty} \P(W) = 0$. For details see Figure \ref{fig:f3}.

\begin{figure}[!hbtp]
\begin{center}
\begin{tabular}{ccc}
\includegraphics[width=5.1cm]{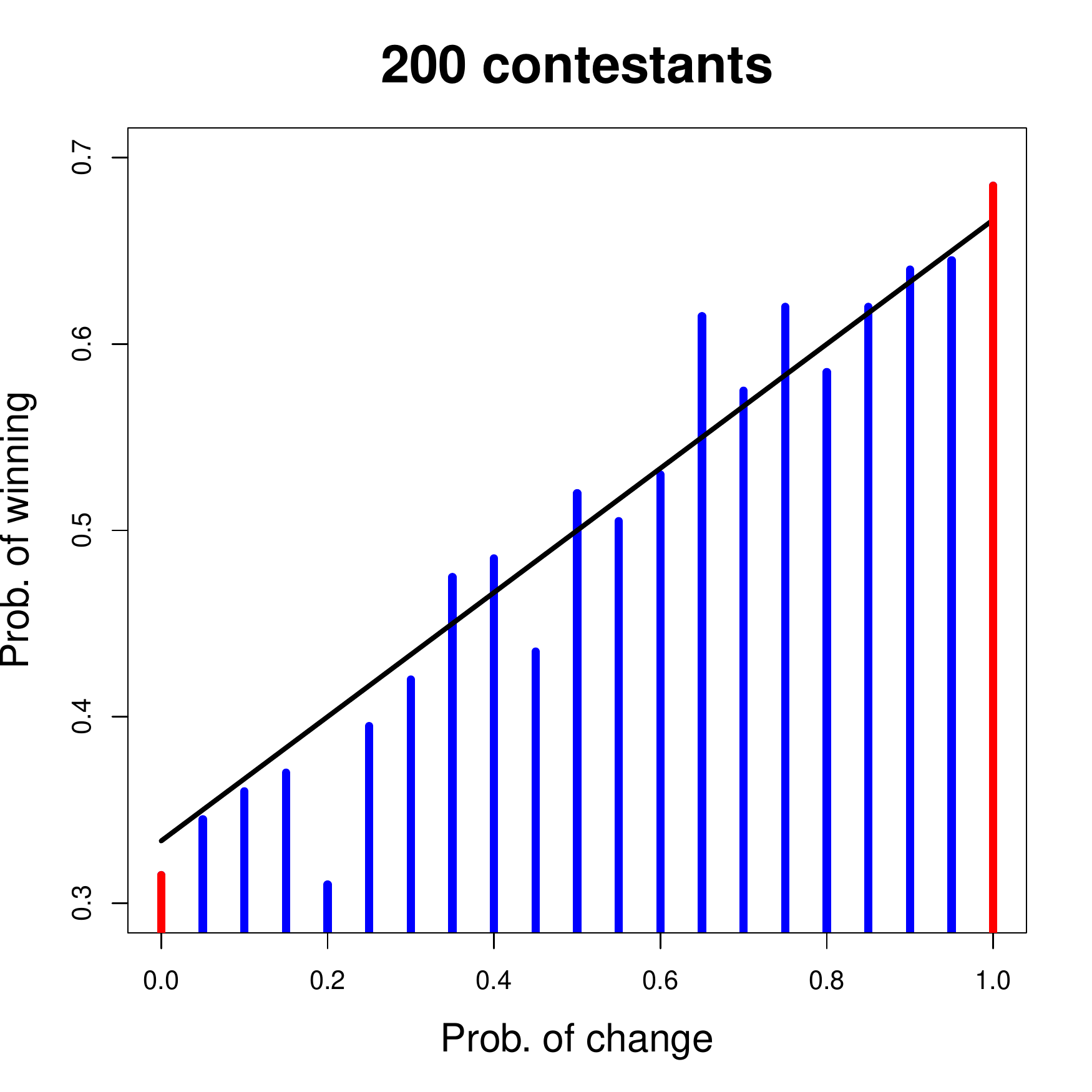} & \includegraphics[width=5.1cm]{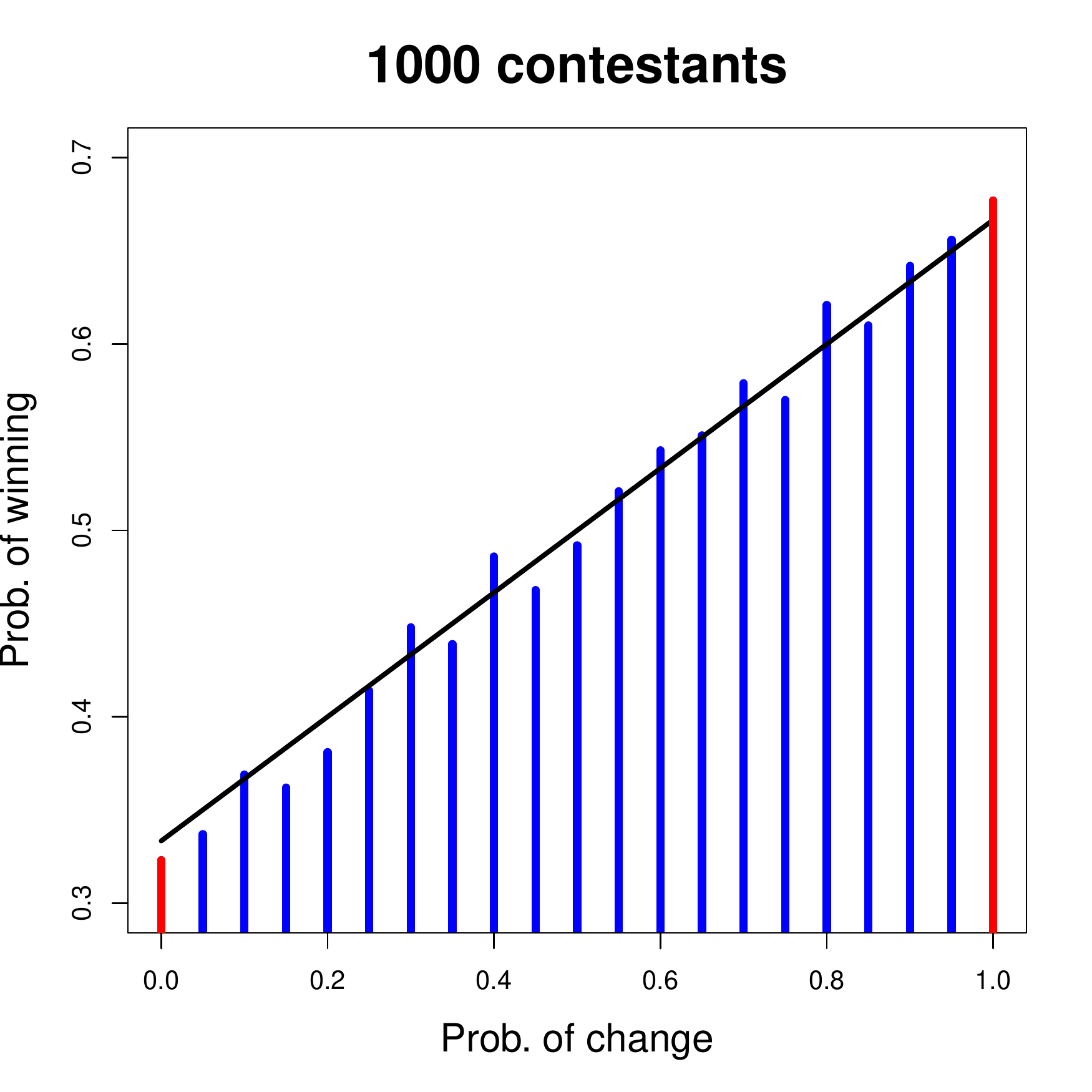} & \includegraphics[width=5.1cm]{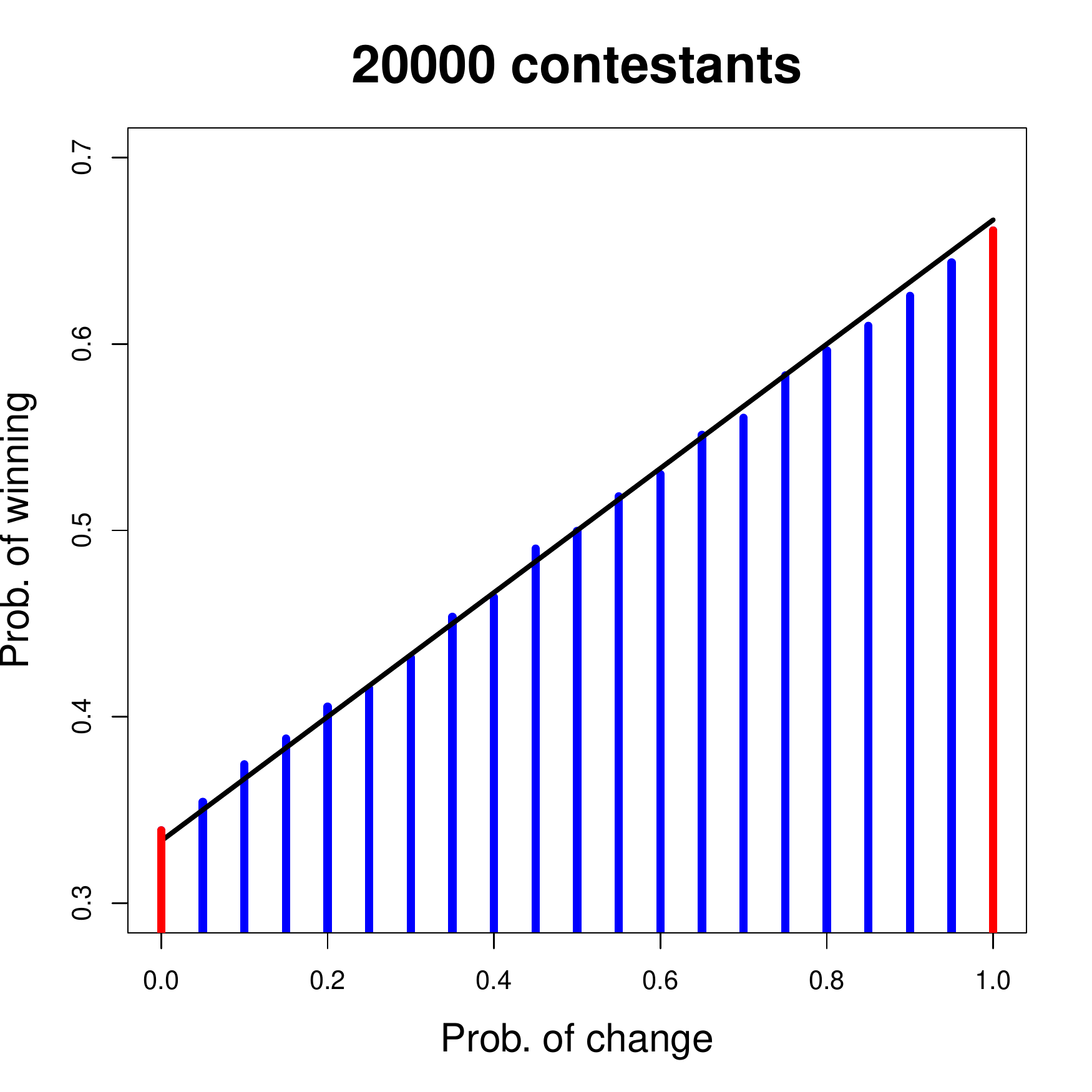}
\end{tabular}
\caption{
Empirical probabilities of winning for the model in Section \ref{sec:MH1}, the straight line is the function of parameter $p$ in equation \eqref{pw1}. Several numbers of steps to check the law of large numbers. } \label{fig:f1}
\end{center}
\end{figure}

\section{Simulation studies}
\label{sec:simu}

\begin{figure}[!hbtp]
\begin{center}
\begin{tabular}{ccc}
\includegraphics[width=5.1cm]{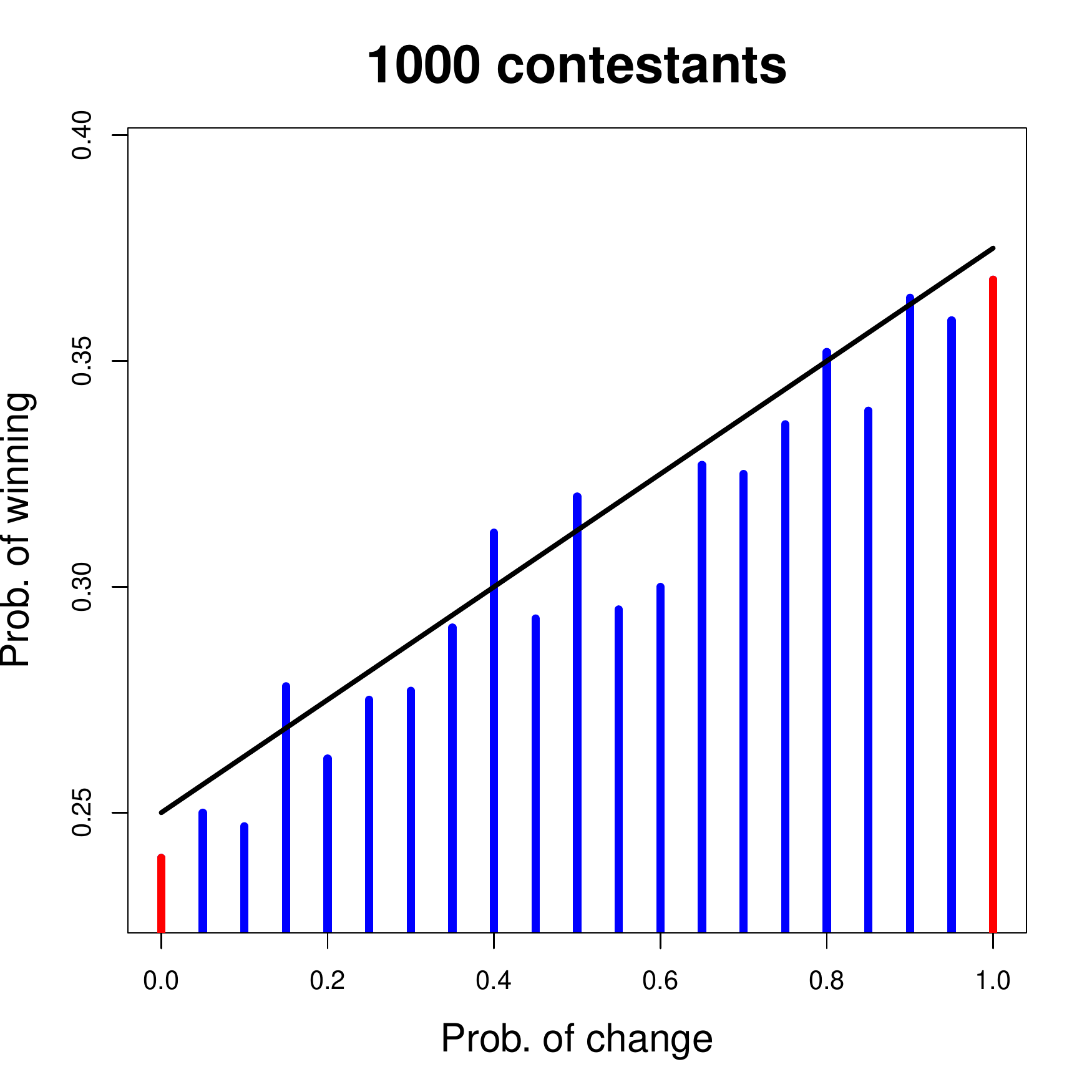} & \includegraphics[width=5.1cm]{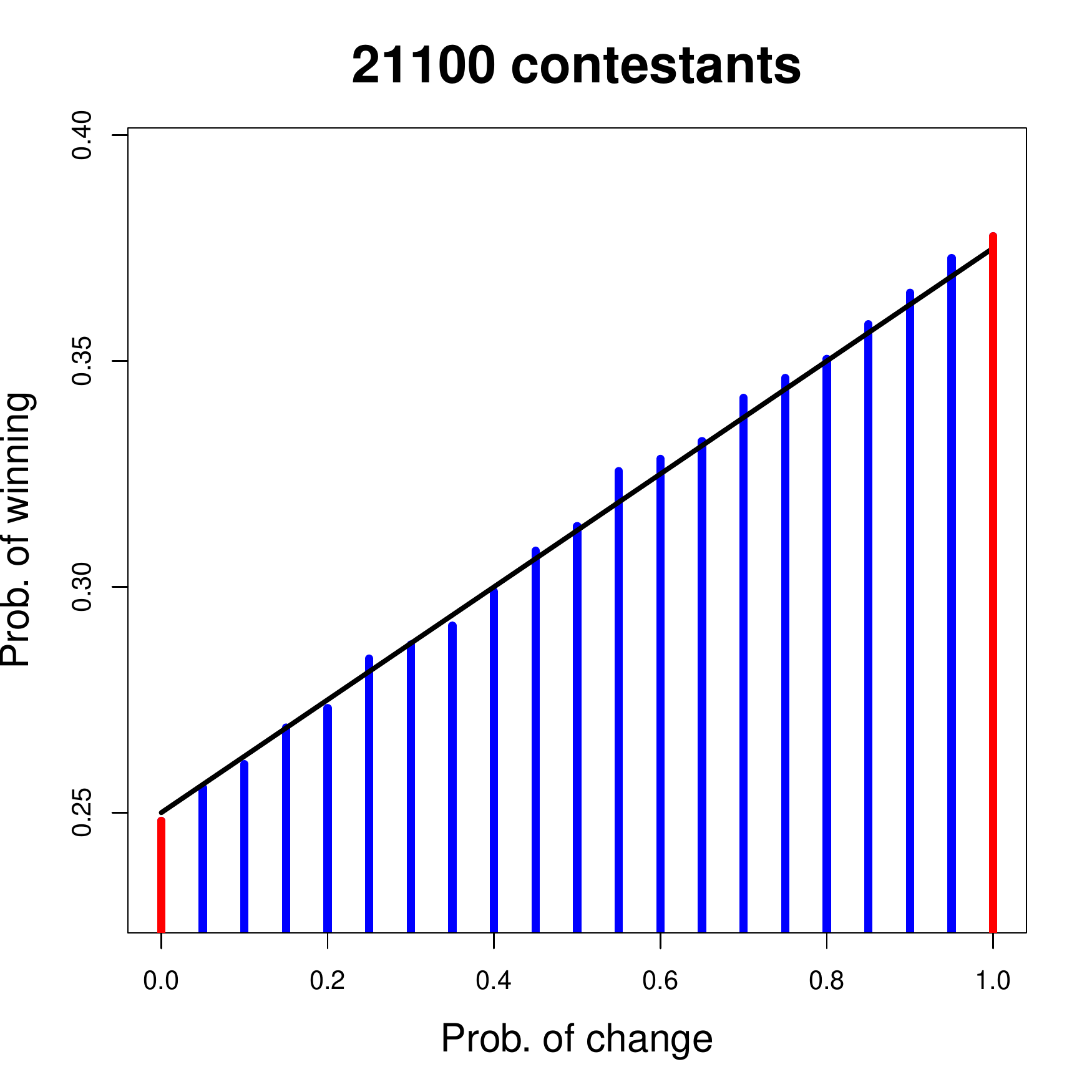} & \includegraphics[width=5.1cm]{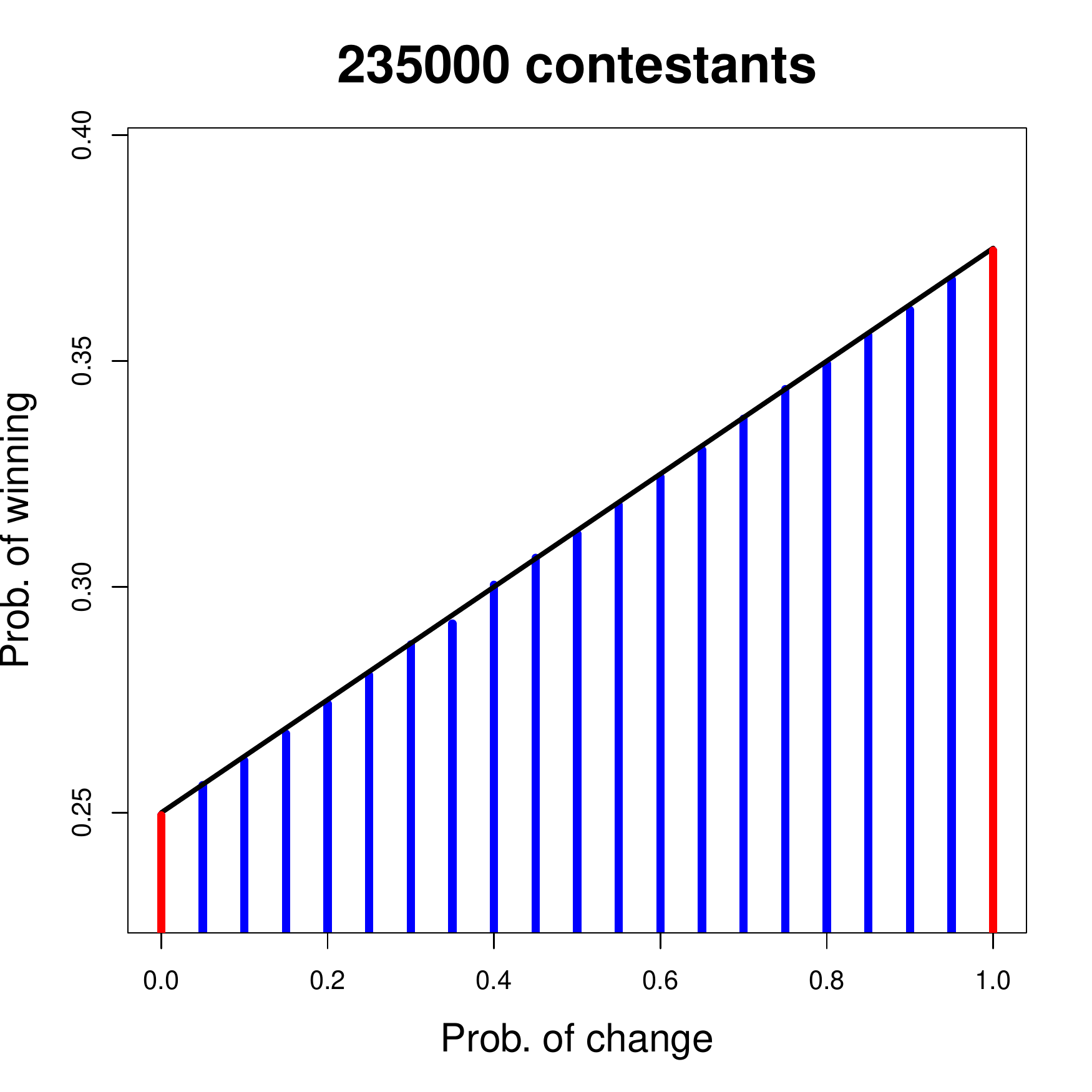}
\\
\includegraphics[width=5.1cm]{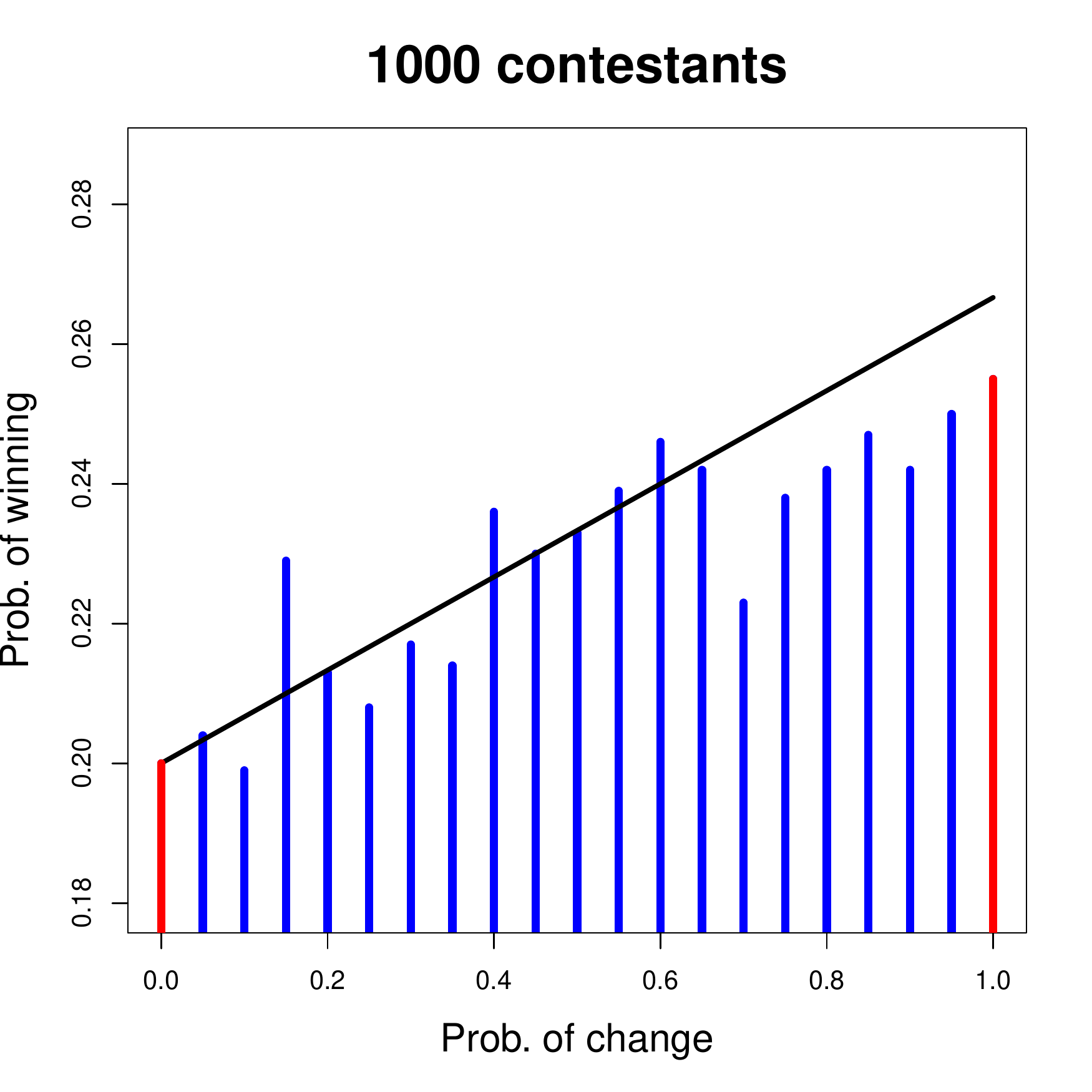} &  \includegraphics[width=5.1cm]{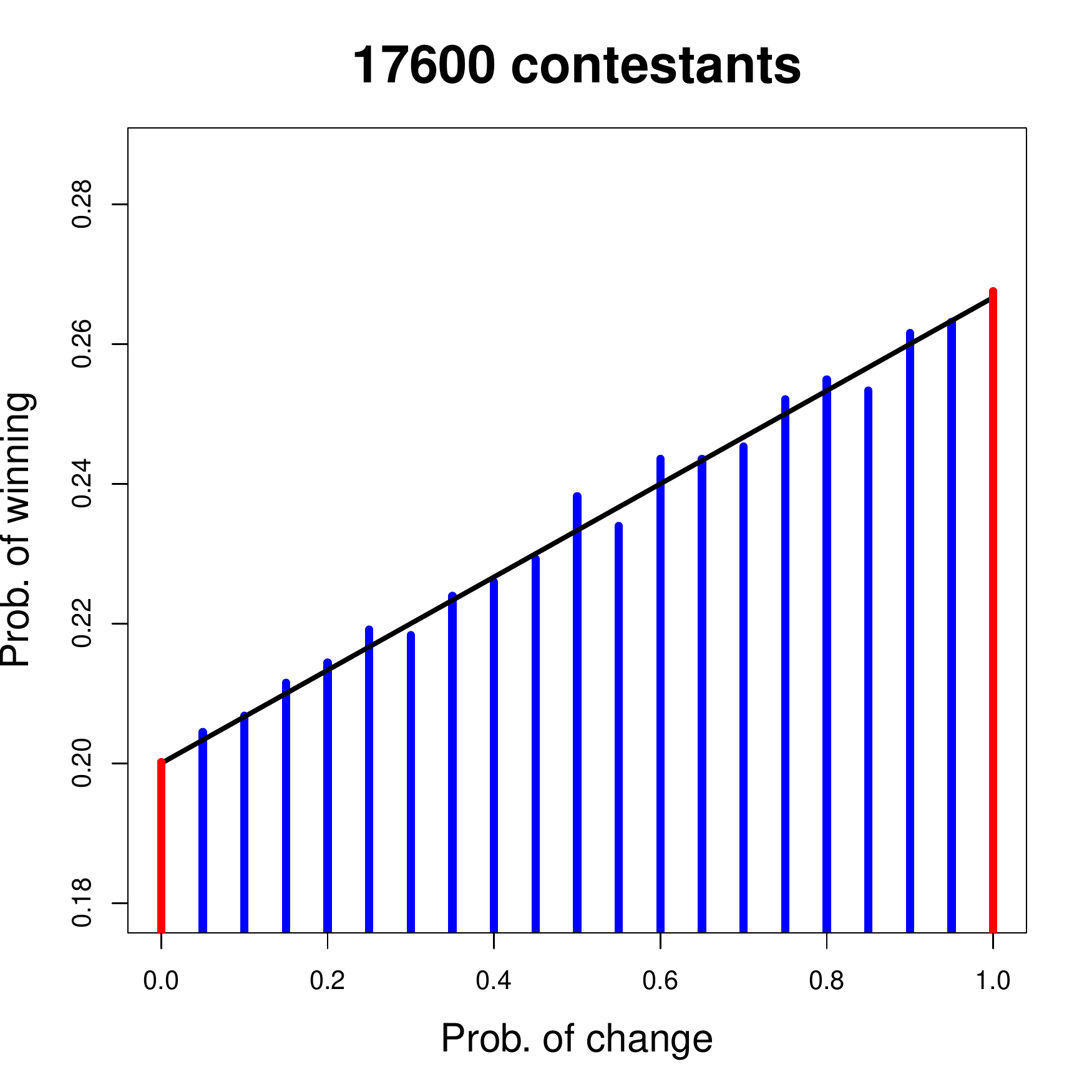} &  \includegraphics[width=5.1cm]{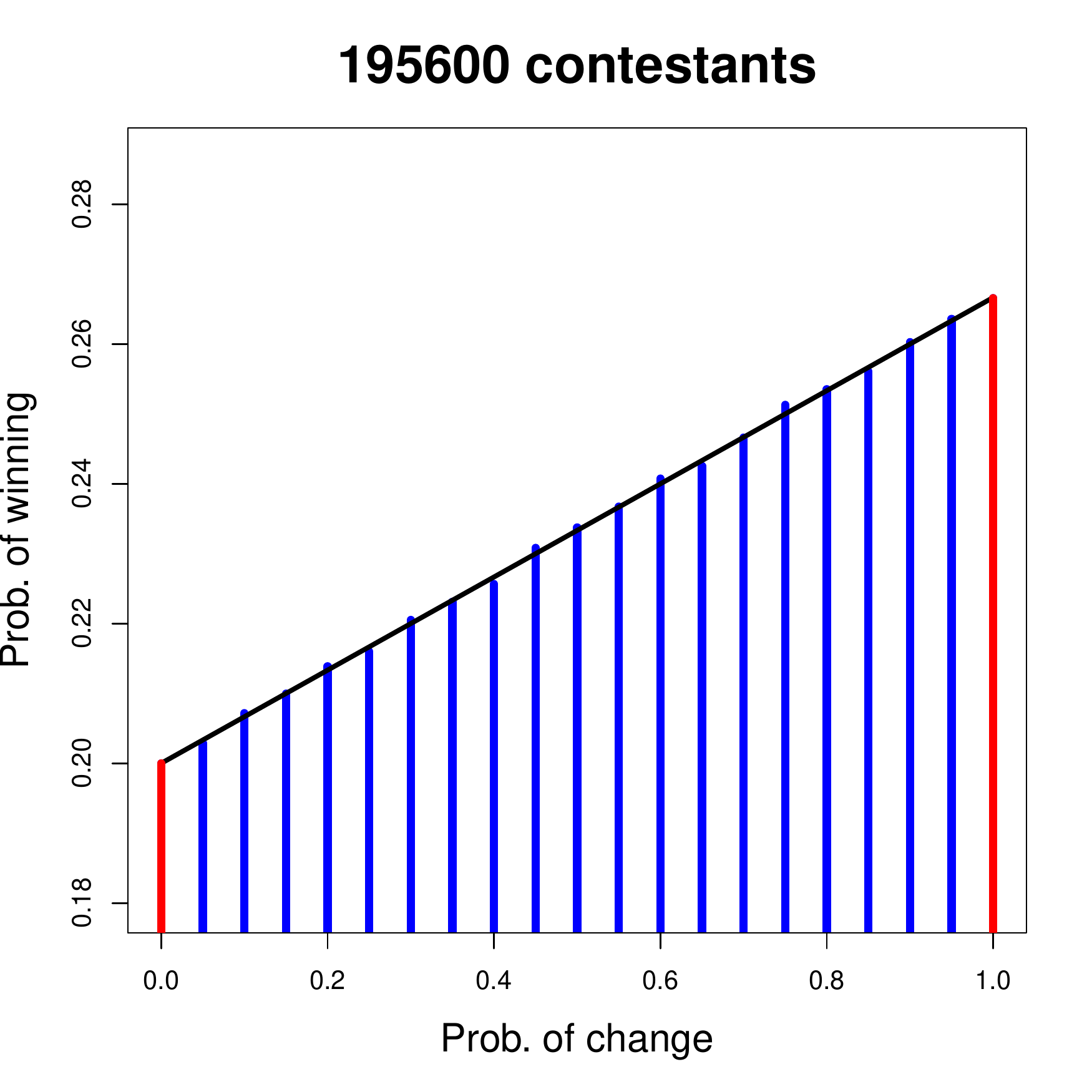}
\\
\includegraphics[width=5.1cm]{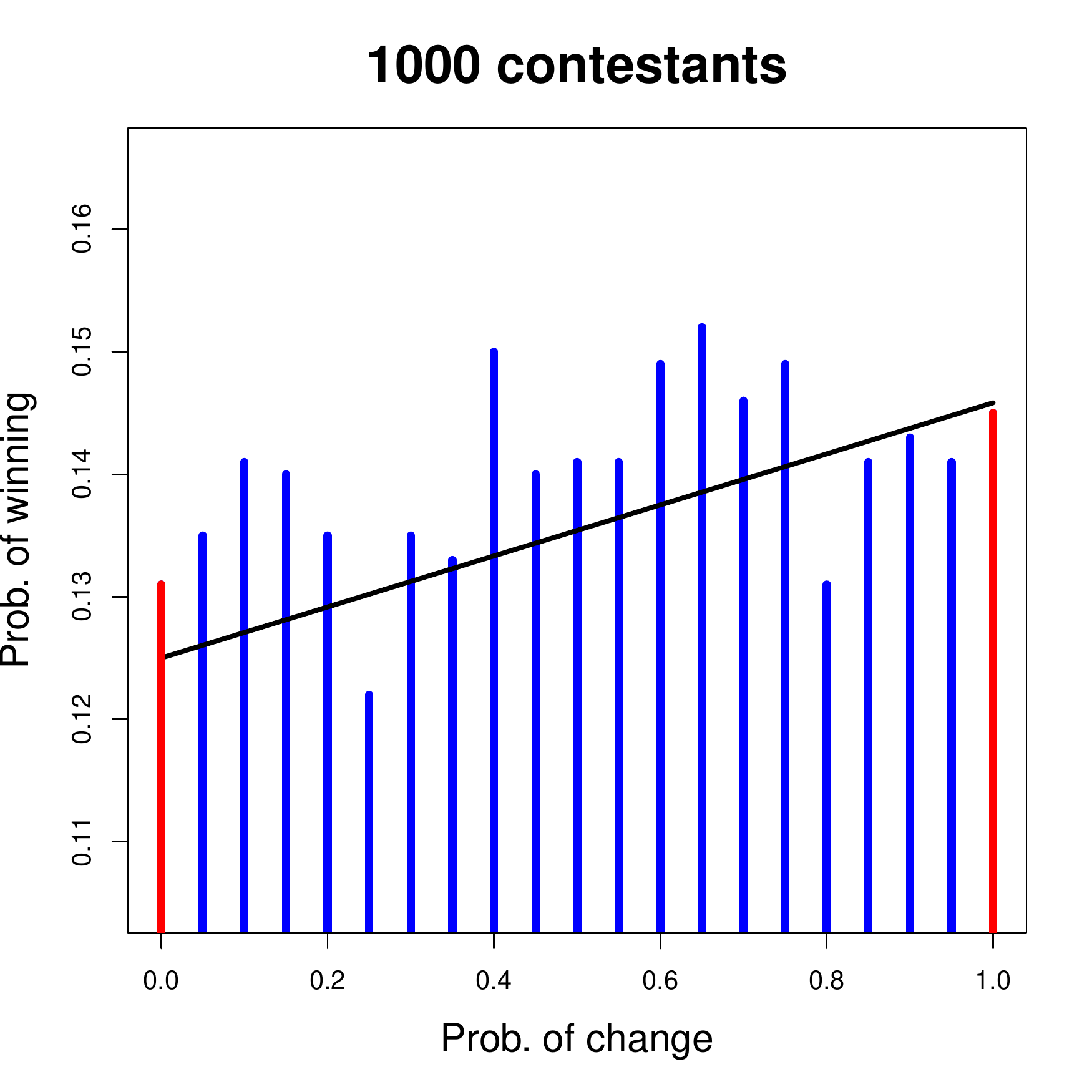} &  \includegraphics[width=5.1cm]{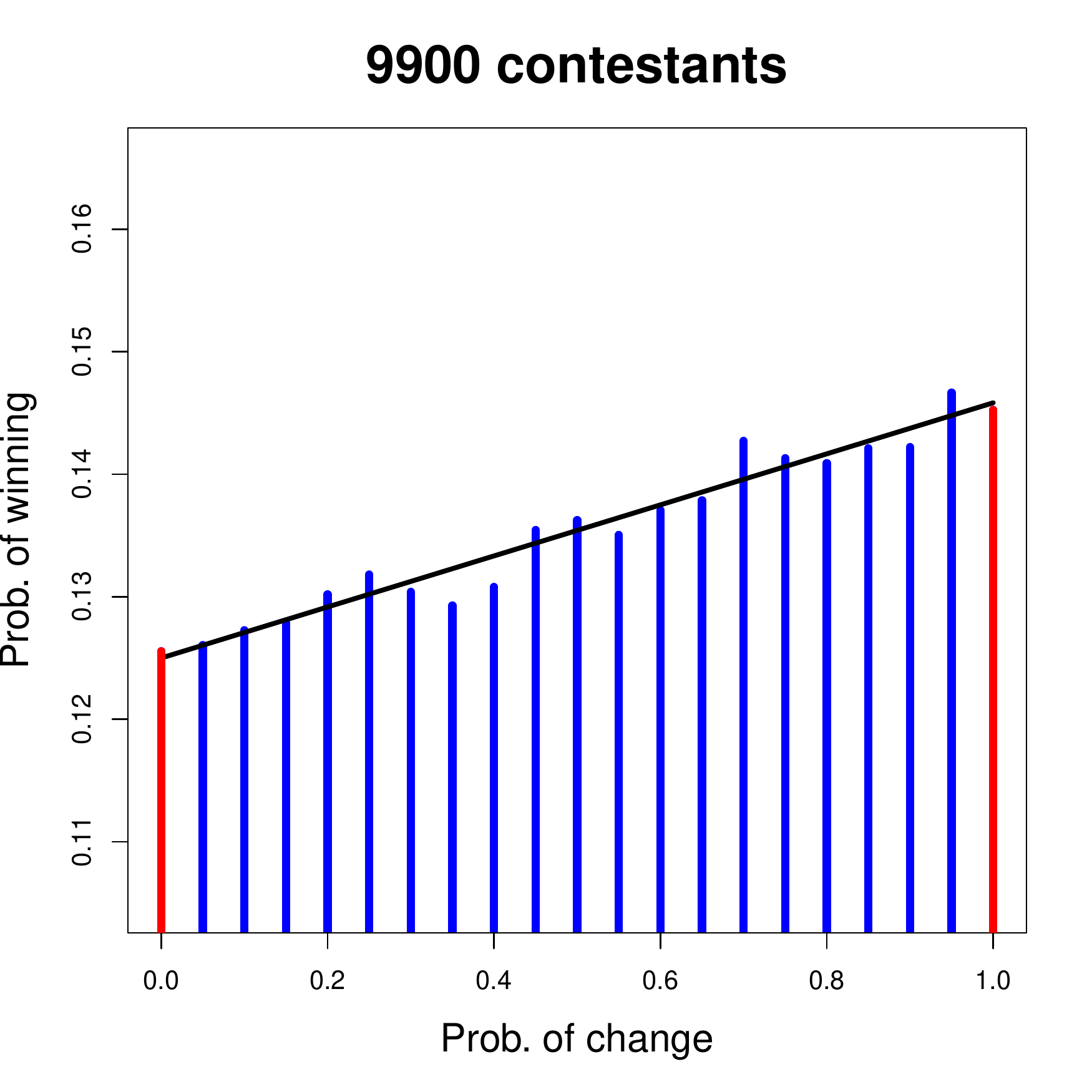} &  \includegraphics[width=5.1cm]{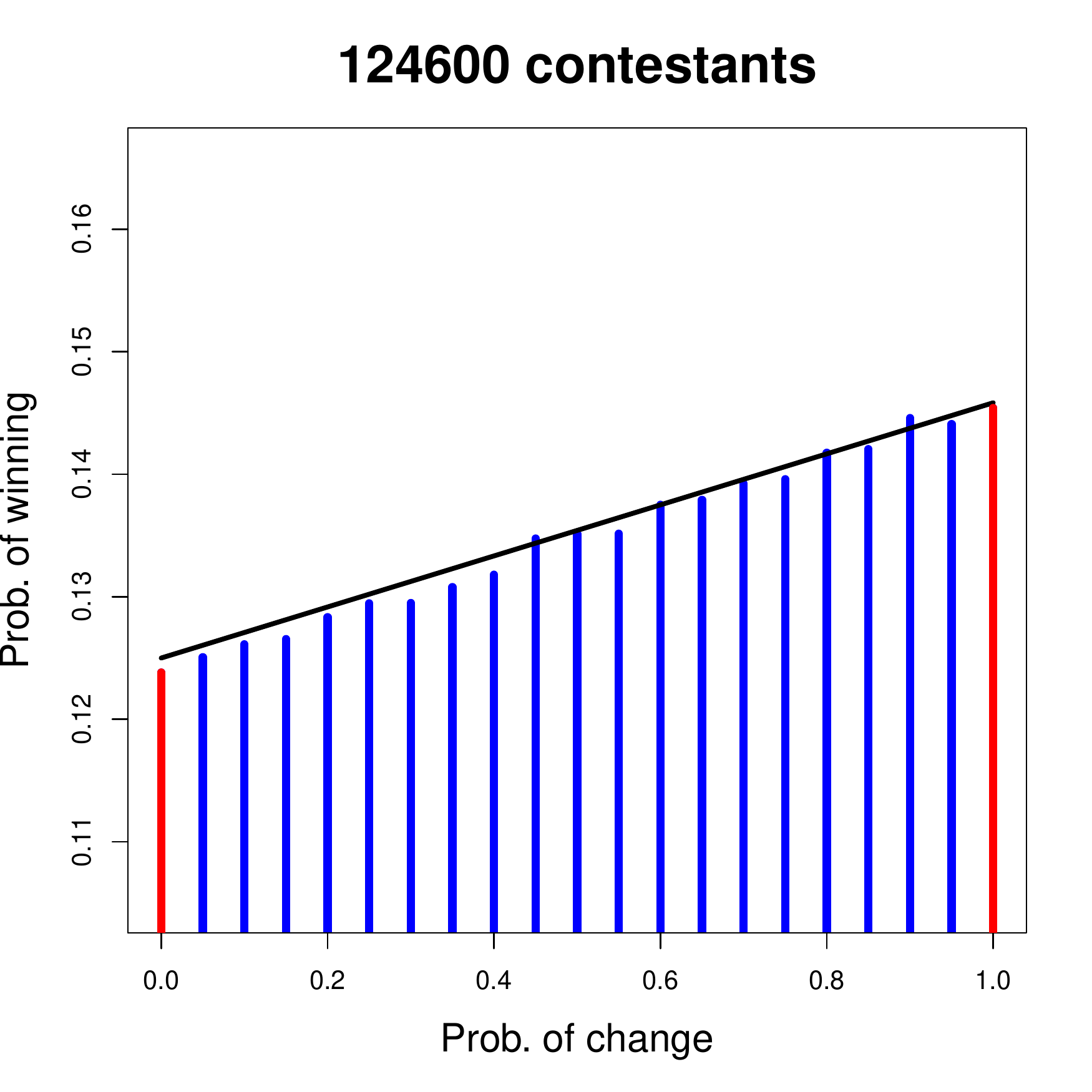}
\\
\includegraphics[width=5.1cm]{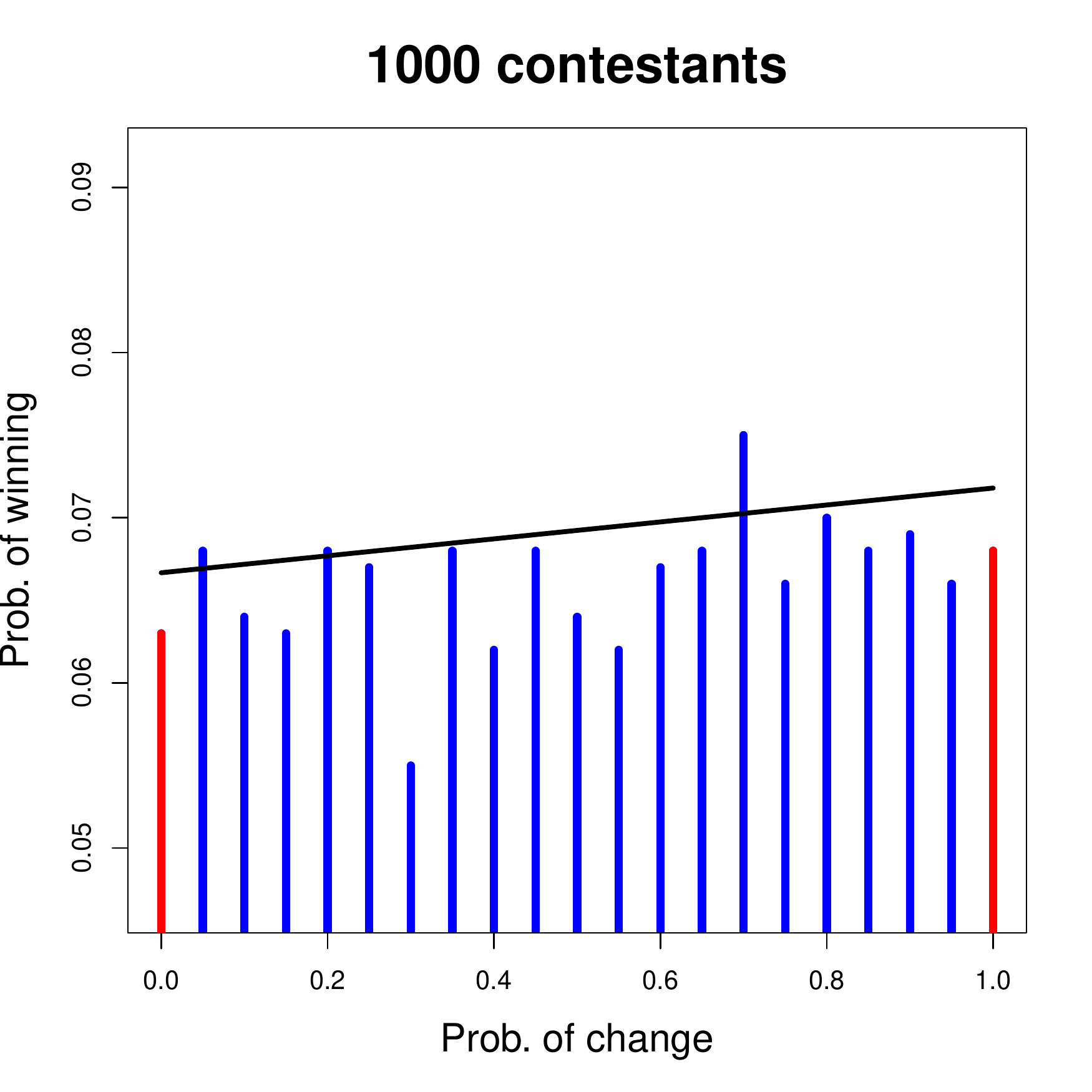} &  \includegraphics[width=5.1cm]{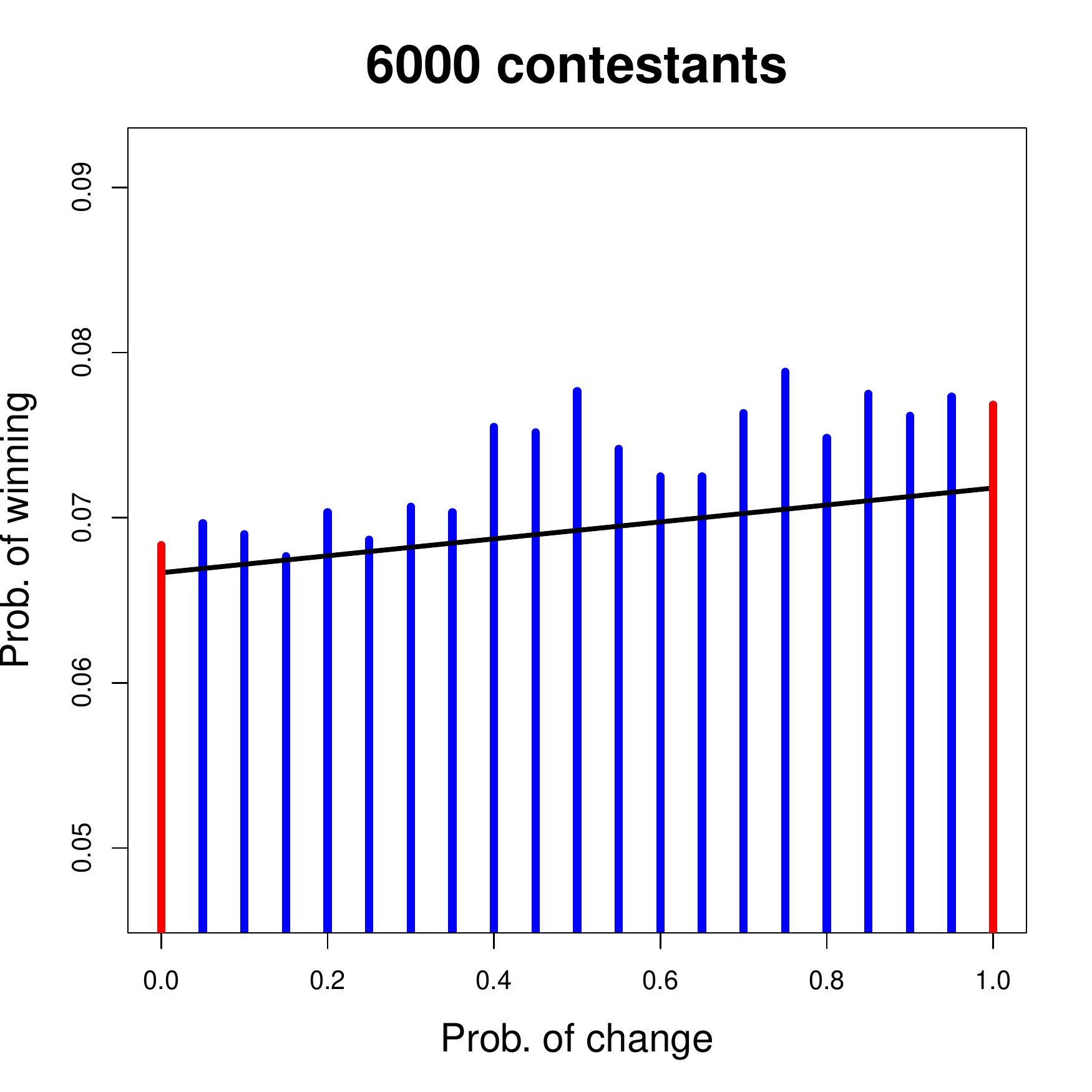} &  \includegraphics[width=5.1cm]{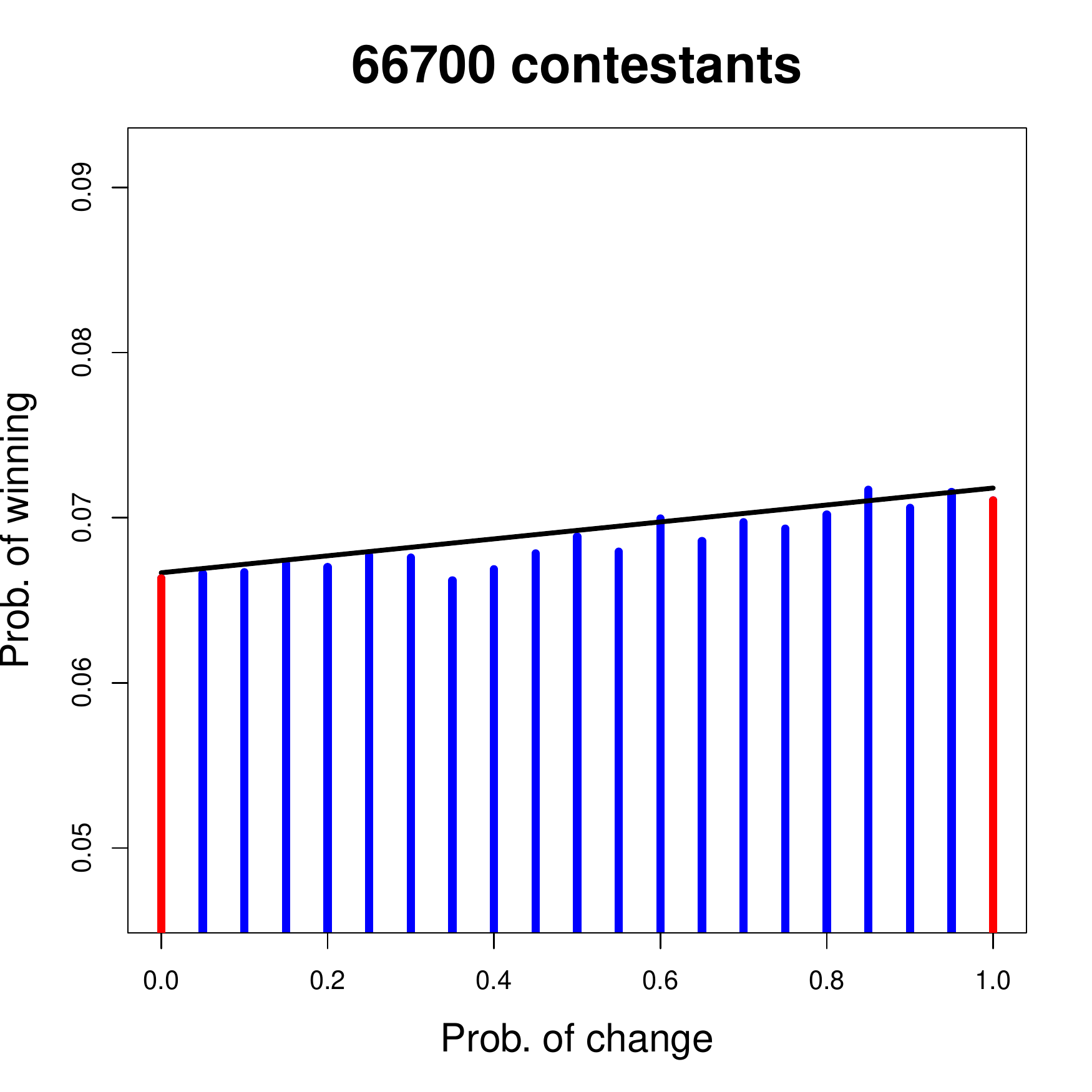}
\end{tabular}
\caption{Generalized Monty Hall problem in Section \ref{sec:MH3}, several number of steps. Second column the bounds by using CLT, last one by using Chebishev inequality. First row with 4 doors, second with 5 doors, third with 8 doors and the last with 15 doors.} \label{fig:f3}
\end{center}
\end{figure}
The main goals of this section are twofold: on the one hand, we analyze the performance of the strategies to check the theoretical solution. On the other hand, we include the arguments to choose the number of iterations that guarantees the empirical probabilities of wining in \eqref{pw1}, \eqref{pw2} and \eqref{pw3} converge to theoretical solutions by the strong law of large numbers (SLLN). The critical values are obtained using the bounds provided by the central limit theorem and Chebyshev inequality.

We develop numerical simulations for the three models in Section \ref{sec:model}. In particular, we analyse the probabilities of winning in \eqref{pw1}, \eqref{pw2} and \eqref{pw3}, as functions of parameter $p$. The empirical probabilities converge to theoretical values by the SLLN. In other words, let a sequence of independent Bernoulli random variables $X_i$, $i= 1, \ldots, l$, where
$$
X_i=
\left\{
\begin{array}{rl}
1,& \text{if the $i$-th contestant win the car} ,\\[0.2cm]
0, & \text{otherwise}.
\end{array}
\right.
$$
By the SLLN we have

$$\displaystyle\frac{1}{l}\displaystyle\sum_{i=1}^l X_i \overset{a.s}{\to} \E(X_i) = \P(W) \ \text{ as } l \to \infty.$$

This means that as long we simulate the process, increasing the number of contestants, the empirical probabilities get closer the theoretical ones. However, we need an additional tool to find a bound for the optimal number of iterations. Then, by using the central limit theorem (CLT) and denoting $\bar{X}_l=\displaystyle\frac{1}{l}\displaystyle\sum_{i=1}^lX_i$, we assume that $ \bar{X}_l\sim N\left(\P(W),\frac{\P(W)(1-\P(W))}{l}\right)$ for $l$ large enough. In other words $\bar{X}_l$  is distributed asymptotically normal with mean $\P(W)$ and variance $\displaystyle\frac{\P(W)(1-\P(W))}{l}$. From the above we obtain the value of $l_0$, such that for all $l > l_0$:

$$ \P\left(\left| \bar{X}_l - \P(W) \right| \geq \varepsilon \right) < \delta.$$

Then, by straightforward calculations we obtain

$$ \P\left( \bar{X}_l  - \P(W) < \varepsilon \right) >1-\frac{\delta}{2}\;\;\Longrightarrow\;\;\P\left(Z< \frac{\varepsilon \sqrt{l}}{\sqrt{\P(W)(1-\P(W))}} \right) >1-\frac{\delta}{2},$$
where $Z \sim N(0,1)$. Denote $x=1-\frac{\delta}{2}$, then find $z_x= \displaystyle\frac{\varepsilon \sqrt{l}}{\sqrt{\P(W)(1-\P(W))}}$. Therefore,
\begin{equation}
\label{lCL}
l^{CL}_0 = \displaystyle\frac{z_x^2\P(W)(1-\P(W))}{\varepsilon^2}.
\end{equation}


We remark that the approximation obtained by the CLT for the probability in \eqref{pw3} is worse if the number of doors $n$ is larger. Then, as an alternative, we use the Chebyshev inequality. Then

\begin{equation}
\label{lCh}
l_0^{Ch} = \frac{\P(W)(1-\P(W))}{\delta \varepsilon^2}
\end{equation}

In each MH problem we use $\varepsilon = \delta = 0.01$ to obtain the corresponding critical $l_0$. For instance, in \eqref{pw1} represented by Figure \ref{fig:f1}, we consider three scenarios with 200, 1000 and 20000 contestants. The last one was obtained using $l_0^{CL}$ in \eqref{lCL}.

The second problem in \eqref{pw2} is represented by Figure \ref{fig:f2}, last column considers the optimal number of iterations using CTL inequality (see \eqref{lCL}). In particular we consider two cases with 5 and 10 doors.

Finally, in \eqref{pw3} represented by Figure \ref{fig:f3}, the second and third scenarios consider the optimal number of iterations using $l_0^{CL}$ and $l_0^{Ch}$, respectively. For this MH version, we consider four cases with 4, 5, 8 and 15 doors.



%
\section{Comments and concluding remarks}
\label{sec:comment}

Some authors have considered the assumption that MH chooses randomly the door with the car (at the beginning). However, this is not relevant. Define
$$A_i : \text { MH puts the car behind door} \ i,$$
and let $\P(A_i ) = \alpha_i$ where $1 \le i \le n$ and $\sum_i \alpha_i = 1$. We are able to prove
$$\P(E) = \sum_1^n \P(E| A_i) \P(A_i) = \sum_1^n \frac{1}{n} \alpha_i = \frac{1}{n}$$
Then Assumptions \ref{assum} are enough.

In addition, note that in the case of Monty Hall problem in Section \ref{sec:MH1}, we obtain $\P(W)=\displaystyle\frac{1}{3} (1+p)$ and in particular, if we assume $p= \displaystyle\frac{1}{2}$, then $\P(W)=\displaystyle\frac{1}{2}$. That is, the solution proposed by the detractors of Vos Savant. In this sense, we claim that the intuition that leads to such solution is the implicit assumption that the contestant does not have a preference to change its initially selected door. In any case, the original question was posed as: Is it better to change or not? Then, Vos Savant solution is absolutely correct.

We remark that for the models in Sections \ref{sec:MH2} and \ref{sec:MH3}, the probabilities of winning respectively \eqref{pw2} and \eqref{pw3}, are linear functions of $p$. However, as $n$ goes to infinity, the slopes of these curves change in different ways. Particularly, \eqref{pw2} converges to $p$, then its slope goes to $1$, and \eqref{pw3} converges to $0$ with slope going to $0$, see Figures \ref{fig:f2}-\ref{fig:f3}.

In a similar sense, we also remark that there exists a critical $n$ such that \eqref{pw3} is strictly smaller than $1/3$, that is the probability \eqref{pw1} when $p=0$.  In this case, the critical value is $n=4$. Then, for all $n > 4 $ the game in Section \ref{sec:MH3} is more advantageous for Monty Hall than the original one (Section \ref{sec:MH1}), see Figure \ref{fig:f3}.

Finally, we highlight the relationship between the original Monty Hall problem and the well-known 3 prisoners paradox (see \cite{car,Gar}). Essentially, the probabilities involved in both problems are the same, where the arguments to calculate them are similar.

\begin{algorithm}{!hbtp}
  \caption{Monty Hall $n$ doors, original version}
  \begin{algorithmic}[1]\label{MHn1}
\STATE fix the number of contestants $l$
\FORALL{ $p \in \mathcal{P}$ in \eqref{Pset}}
\FORALL{ $i = 1, \ldots, l$}
\STATE fix the door with the car (for instance $1$)
\STATE let $W=0$ (initial number of winnings)
\STATE choose $X \in \{1, \ldots ,n\}$ uniformly (initial selected door)
\IF{ $X=1$ }
\STATE choose $Y \in \{2, \ldots , n\}$ uniformly
\ELSE
\STATE choose $Y = 1$
\ENDIF
\STATE let $(X,Y)$ closed doors
\STATE sample $C \sim Bern(p)$
\IF{ $C=0$ }
\STATE $Z=X$ (final selected door)
\ELSE
\STATE $Z=Y$
\ENDIF
\IF{ $Z=1$ }
\STATE $W=W+1$
\ELSE
\STATE $W=W$
\ENDIF
\ENDFOR
\ENDFOR
\STATE plot values of arrays $\displaystyle\left(p, \frac{W}{l}\right)$
\end{algorithmic}
\end{algorithm}









\section{Appendix}
\label{sec:app}

Finally, we include the pseudo-codes. For that, we need to define

\begin{equation}
\label{Pset}
\mathcal{P} = \{ 0.05\cdot k: 0 \le k \le 20\}.
\end{equation}

Then, the Algorithm \ref{MHn1} is used for models in Sections \ref{sec:MH1} and \ref{sec:MH2}. The Algorithm \ref{MHn2} is related to the model in Section \ref{sec:MH3}.

\bigskip

\bigskip
\begin{algorithm}{!hbtp}
  \caption{Monty Hall $n$ doors, opening one}
  \begin{algorithmic}[1]\label{MHn2}
\STATE fix the number of contestants $l$
\FORALL{ $p \in \mathcal{P}$ in \eqref{Pset}}
\FORALL{ $i = 1, \ldots, l$}
\STATE fix the door with the car (for instance $1$)
\STATE let $W=0$ (initial number of winnings)
\STATE choose $X \in \{1, \ldots ,n\}$ uniformly (initial selected door)
\IF{ $X=1$ }
\STATE choose $Y \in \{2, \ldots , n\}$ uniformly
\ELSE
\STATE choose $Y \in \{2, \ldots, X-1, X+1, \ldots , n\}$ uniformly
\ENDIF
\STATE let $\mathcal{R}=\{1, \ldots ,n\} \setminus \{X,Y\}$ (remaining closed doors)
\STATE sample $C \sim Bern(p)$
\IF{ $C=0$ }
\STATE $Z=X$ (final selected door)
\ELSE
\STATE choose $Z\in \mathcal{R}$ uniformly
\ENDIF
\IF{ $Z=1$ }
\STATE $W=W+1$
\ELSE
\STATE $W=W$
\ENDIF
\ENDFOR
\ENDFOR
\STATE plot values of arrays $\displaystyle\left(p, \frac{W}{l}\right)$
\end{algorithmic}
\end{algorithm}

\subsection*{Acknowledgements}
CCC was partially supported by DIUBB 173408 2/I of the Universidad del B\'io-B\'io. MGN is supported by Fondo Especial DIUBB 1901083-RS.

\end{document}